\documentclass[reqno,11pt]{amsart}
\usepackage[dvipsnames]{xcolor}

\usepackage{amsmath}
\usepackage{amsfonts}
\usepackage{amssymb}
\usepackage{amsthm}
\usepackage{amsxtra}
\usepackage{mathtools}
\mathtoolsset{centercolon}		% damit schaut := besser aus ohne \coloneqq zu verwenden 
\usepackage{mathrsfs}			  % caligraphic symbols
\usepackage{empheq}
\numberwithin{equation}{section}

\usepackage{enumitem}			  % customization of enumerations and itemize
\usepackage{esint}						% for \fint - struck through integral (average)

\usepackage{graphicx}
\usepackage{tikz}
\usetikzlibrary{hobby, intersections, calc}
\usetikzlibrary{decorations.pathmorphing}

\newtheorem{theorem}{Theorem}[section]

\newtheorem{lemma}[theorem]{Lemma}
\newtheorem{lem}[theorem]{Lemma}
\newtheorem{proposition}[theorem]{Proposition}
\newtheorem{corollary}[theorem]{Corollary}
\newtheorem{remark}[theorem]{Remark}

\newtheorem{assumption}{Assumption}[section]

\newcommand{\name}[1]{#1}			% for names of authors
\newcommand{\df}[1]{\textit{#1}}				  % for definitions 

\AtBeginDocument{
	\let\div\relax 
	\DeclareMathOperator{\div}{div}

}

\newcommand{\R}{\mathbb{R}}
\newcommand{\N}{\mathbb{N}}

\newcommand{\C}{\mathbb{C}}
\renewcommand{\S}{\mathcal{S}}

\newcommand{\Rigid}{\mathcal{R}}						  % Set of rigid motions
\newcommand*\dd{\mathop{}\!\mathrm{d}}

\newcommand{\pv}{\textup{pv-}\hspace{-4pt}}		% Chauch principle value

\newcommand{\blank}{{\,\cdot\,}}

\newcommand{\epsi}{\varepsilon}
\newcommand{\eps}{\varepsilon}

\newcommand{\D}{\nabla}
\newcommand{\Hd}{\mathcal{H}}				% Hausdorff measure
\newcommand{\wto}{\rightharpoonup} 		%weak to

\newcommand{\dxd}{ {d \times d} }
\newcommand{\dxdxd}{ {d \times d \times d} }

\DeclareMathOperator{\dist}{dist}
\DeclareMathOperator{\supp}{supp}

\DeclareMathOperator{\tr}{tr}

\newcommand{\Dd}{\mathbb{D}} 		 % for the metric 
\newcommand{\Dsc}{\mathscr{D}} 		% script D for the nonlocal displacement
\newcommand{\Dfr}{\mathfrak{D}}		 % fraktur D for the nonlocal divergence

\newcommand{\Sn}{\mathcal{S}_{n}}
\newcommand{\OO}{{\widetilde{\Omega}}}

\usepackage{hyperref}
\hypersetup{
	colorlinks=true, 
	citecolor = Green}

\title{Nonlocal-to-local limit in  linearized  viscoelasticity} 
\date{\today}

\begin{document}

	%BEGIN_FOLD Authors and stuff
	
	\author[M. Friedrich]{Manuel Friedrich} 
	\address[Manuel Friedrich]{Department of Mathematics, Friedrich-Alexander Universit\"at Erlangen-N\"urnberg. Cauerstr.~11,
		D-91058 Erlangen, Germany, \& Mathematics M\"{u}nster,  
		University of M\"{u}nster, Einsteinstr.~62, D-48149 M\"{u}nster, Germany.
	}	
	\email{manuel.friedrich@fau.de}
	\urladdr{https://www.math.fau.de/angewandte-mathematik-1/mitarbeiter/prof-dr-manuel-friedrich/}

	\author[M. Seitz]{Manuel Seitz} 
	\address[Manuel Seitz]{Vienna School of Mathematics, University of
		Vienna, Oskar-Morgenstern-Platz~1, A-1090 Vienna, Austria, \& Faculty of Mathematics,  University of
		Vienna, Oskar-Morgenstern-Platz~1, A-1090 Vienna, Austria.
	}
	\email{manuel.seitz@univie.ac.at}

	\author[U. Stefanelli]{Ulisse Stefanelli} 
	\address[Ulisse Stefanelli]{Faculty of Mathematics, University of
		Vienna, Oskar-Morgenstern-Platz 1, A-1090 Vienna, Austria,
		Vienna Research Platform on Accelerating
		Photoreaction Discovery, University of Vienna, W\"ahringerstra\ss e 17, 1090 Wien, Austria,
		\& Istituto di
		Matematica Applicata e Tecnologie Informatiche {\it E. Magenes}, via
		Ferrata 1, I-27100 Pavia, Italy.
	}
	\email{ulisse.stefanelli@univie.ac.at}
	\urladdr{http://www.mat.univie.ac.at/$\sim$stefanelli}
	
	\keywords{Peridynamics, viscoelasticity, Kelvin-Voigt rheology, nonlocal-to-local limit, evolutionary $\Gamma$-convergence.}

	\begin{abstract}
		We study the quasistatic evolution of a linear peridynamic Kelvin-Voigt viscoelastic material. More specifically, we consider  the gradient flow of a nonlocal elastic energy with respect to a nonlocal viscous dissipation.  Following an evolutionary $\Gamma$-convergence approach, we prove that the solutions  of  the nonlocal problem converge to the solution of the   local     problem, when the peridynamic horizon tends to $0$, that is, in the nonlocal-to-local limit. 
	\end{abstract}

	\subjclass[2020]{
		74A70,  % peridynamics
		70G75, % Variational methods for problems in mechanics
		35Q74, % PDEs in connection with mechanics of deformable solids
		49J45, %Methods involving semicontinuity and convergence; relaxation
		74D05. %Linear constitutive equations for materials with memory
		}   
	\maketitle

	%%%%%%%%%%%%%%%%%%%%%%%%%%%%%%%%%%%%%%%%%%%%%%%%%%%%
	\section{Introduction}
	%%%%%%%%%%%%%%%%%%%%%%%%%%%%%%%%%%%%%%%%%%%%%%%%%%%%%

   	\noindent\textbf{The model.}   	In this paper, we are interested in the quasistatic evolution   of linear     peridynamic \name{Kelvin}--\name{Voigt} viscoelastic materials.   Let      $\Omega \subset \R^d$ ($d \geq 2$) be a bounded reference domain and   indicate by  $u\colon [0,T] \times \Omega \to \R^d$ the     \df{time dependent displacement}. To prescribe nonlocal homogeneous Dirichlet boundary values, we fix  $\OO \in \R^d$ such that $\Omega \subset \subset \OO$ (compactly contained) and extend  $u$ by setting $u(x) = 0$  for $x \in \OO \setminus \Omega$. Following \cite{mengeshaVariationalLimitClass2015} and moving within the framework of (state-based) peridynamics, we define the \df{nonlocal    elastic energy}    of     a homogeneous, isotropic material as 
	\begin{align}\label{eq:energy}
		E_{\rho}(u) = &\frac{\beta}{2} \int_\OO (\Dfr_\rho(u)(x))^2 \dd x \nonumber \\
		+ &\frac{\alpha}{2} \int_\OO \int_\OO \rho(x'-x) \left(\Dsc(u)(x',x) -   \frac{1}{d}     \, \Dfr_\rho(u)(x)\right)^2 \dd x' \dd x.  
	\end{align}
	Here, $\rho : \R^d \to [0, \infty)$ is a radially symmetric kernel, weighting the strength of nonlocal    interactions.    
	We assume that $\rho (\xi ) = 0$ for $\vert \xi \vert >   h    $, where   $h$     is the so-called \df{peridynamic horizon}. Furthermore,  we let   
	\begin{equation}\label{eq:nonlocal strain}
		\Dsc(u)(x',x) := \frac{(u(x') - u(x)) \cdot (x'-x)}{| x' - x|^2}
	\end{equation}
	 indicate   the \df{nonlocal strain}   in    the direction $(x'-x)/| x' - x|$.   The   weighted nonlocal strain, called \df{nonlocal divergence}, is given by
	\begin{equation}\label{eq:nonlocal divergence}
		\Dfr_\rho(u)(x) = \pv\int_\OO \rho(x'-x) \Dsc(u)(x',x) \dd x',
	\end{equation}
	where  'pv' indicates  the Cauchy principle value.  In fact,  it can be shown that   $\Dfr_\rho(u) \to \div u$     as $\rho \to \delta_0$  (Dirac  delta  in $0$) suitably,  see Lemma~\ref{lem:convergence of strain} below. 	The constants $\alpha, \beta > 0$ are material parameters related to the shear and bulk modulus, respectively. 	
	
	  Under   suitable assumptions on the kernel $\rho$,        \name{Mengesha} and \name{Du} showed that $E_{\rho}$ admits a unique minimizer \cite[Theorem~1.1]{mengeshaVariationalLimitClass2015}. Moreover, they proved that $E_\rho$ $\Gamma$-converges \cite{DalMaso_Gamma} with respect to the strong $L^2$-topology to the \df{local elastic energy}
	\begin{align}\label{eq:gamma limit of energy}
		E(u) = \frac{\lambda}{2} \int_\Omega \div u(x)^2 \dd x + \frac{\mu}{2} \int_\Omega \vert \epsi(u(x)) \vert^2 \dd x
	\end{align}
	when passing from the nonlocal to the local problem, that is, when  $\rho \to \delta_0$ suitably   \cite[Theorem~1.2]{mengeshaVariationalLimitClass2015}.  Here $\epsi(u) = (\D u + \D u^\top)/2$ denotes the \df{symmetric gradient of $u$}.   The \df{Lam\'e coefficients} $\lambda$ and $\mu$ can be expressed in terms of $\alpha$, $\beta$ and the dimension $d$,  see \eqref{eq:tensors intro} below  or  \cite[Appendix~A]{mengeshaVariationalLimitClass2015}.  Note that the integration domain in \eqref{eq:gamma limit of energy} is merely $\Omega$, for admissible configurations for $E$ will be considered as belonging to $H^1_0(\Omega;\R^d)$.

 	In order to incorporate viscous dissipation in the model,  we introduce the \df{nonlocal potential of   viscous      forces} given by
	\begin{align}\label{eq:dissipation potential}
		D_\rho(\dot u) = \frac12 \int_\OO \int_\OO \rho(x'- x) \vert \Dsc(\dot u) (x',x)  \vert^2 \dd x' \dd x,
	\end{align}
	where, here and in the following,  the symbol  $\dot u$ denotes the time-derivative of $u$.

	As it turns out, the energy $E_\rho$ is finite if and only if 
	\begin{equation*}%\label{eq:def nonlocal seminorm intro}
		\vert u  \vert_{\rho}^2 :=  \int_{\OO} \int_{\OO} \rho(x'-x) \vert \Dsc(u) (x',x) \vert^2 \dd x' \dd x 
	\end{equation*}
	is finite. This motivates the definition of the following  \df{nonlocal energy space} 
	\begin{align*}
		  \S_{\rho} (\OO;\R^d) := \left\{ u \in L^2(\OO; \R^d) : \vert u \vert_{\rho} < \infty \textup{  and  } u = 0  \textup{ on }\OO\setminus \Omega  \right\}.     
	\end{align*}
	 Due to the nonlocal boundary values, a    Poincar\'e inequality is available, and thus    $\S_{\rho}(\OO;\R^d)$ is a Hilbert space when endowed with the inner product 
	\begin{align*}
		( u, v )_{\rho}  := \int_\OO \int_\OO \rho(x' -x ) \Dsc(u)(x',x) \Dsc(v)(x',x) \dd x' \dd x.
	\end{align*}
	
\emergencystretch3em
	  \noindent\textbf{Our result.}   Given some \df{initial displacement} $u_{\rho,0} \in \S_\rho(\OO;\R^d)$, we look for a solution $u_\rho \in   H^1   (0,T;\S_{\rho}(\OO;\R^d))$ to the \df{nonlocal viscoelastic problem} given by
	\begin{subequations}	  \label{eq:nonlocal problem intro}
		\begin{empheq}[left=\empheqlbrace]{alignat = 2}
		 	\dd_\rho   D_\rho (\dot u_\rho) + \dd_\rho  E_\rho (u_\rho) &= 0   &\quad &\text{in  } \S_{\rho}(\OO;\R^d)^* \quad     \textup{ for  almost every } t \in [0,T],        \label{eq:nonlocal problem intro 1} \\
			u_\rho(0,\cdot) &= u_{\rho,0}  &\quad &\text{in }  \OO.
		\end{empheq}  
	\end{subequations}
  Here, $\dd_\rho D_\rho, \, \dd_\rho E_\rho \colon \S_\rho(\OO; \R^d) \to \S_\rho(\OO; \R^d)^*$ (dual space) are the Fr\'echet differentials of the quadratic functionals $D_\rho$ and $E_\rho$, see Corollary~\ref{cor:frechet}.

	In particular, $\dd_\rho D_\rho$ and $\dd_\rho E_\rho$ are linear, and existence and uniqueness of solutions to the nonlocal   viscoelastic     problem follow from classical theory,    see, e.g.,  \cite[Section~6.4]{kruzikMathematicalMethodsContinuum2019} or \cite[Section~III.~2]{showalterMonotoneOperatorsBanach2013}.   %  \cite[Theorem~IV~6.1]{showalterMonotoneOperatorsBanach2013}.  

	The main  goal of the paper is to show that the solutions to the nonlocal   viscoelastic   problem converge to a solution of the \df{local   viscoelastic   problem} in the nonlocal-to-local limit, that is, when the    kernels $\rho$ converge  to  $\delta_0$      in a suitable sense. 
	The \df{local    viscoelastic   problem} is  that of finding   a weak solution       $u \in    H^1    (0,T;H^1_0(\Omega;\R^d))$ of
	\begin{subequations}\label{eq:pde intro}
		\begin{empheq}[left=\empheqlbrace]{alignat = 2}
			-\D \cdot (\Dd \epsi(\dot u) + \C \epsi(u)) &= 0 &\qquad  &\textup{a.e.\ in  } (0,T)\times \Omega, \label{eq:pde intro1} \\
			u(0,\cdot) &= u_0  &\qquad &\textup{in  } \Omega,
		\end{empheq}  
	\end{subequations}
	where  $u_0 \in H^1_0(\Omega;\R^d)$  is a given initial displacement.  
	  Here,   $\C$   and $\Dd$  are  the \df{elasticity tensor}  and the  \df{viscosity tensor}, respectively,  given by
	\begin{align}
		\C_{ijkl} &=   \frac{2\alpha}{d+2}     \delta_{ik} \delta_{jl} +   \left( \beta -  \frac{2\alpha}{d(d+2)}  \right)     \delta_{ij} \delta_{kl} = \mu  \delta_{ik} \delta_{jl} + \lambda  \delta_{ij} \delta_{kl} , \label{eq:tensors intro} \\
		\Dd_{ijkl} &=   \frac{2}{d+2}     \delta_{ik} \delta_{jl} +   \frac{1}{d+2}      \delta_{ij} \delta_{kl}, \label{eq:tensor intro D}
	\end{align}
	 for $i,j,k,l  = 1, \dots, d$,   where     $\delta_{ij}$     is          the Kronecker delta.  Note that   we can rewrite  the local elastic  energy \eqref{eq:gamma limit of energy} as
	\begin{align}\label{eq:loc energy intro}
		E(u) = \frac12 \int_\Omega \C \epsi(u(x)) : \epsi(u(x)) \dd x, \qquad 
	\end{align}
	see  Lemma~\ref{lem:energy as tensor}   below.  Moreover, we set 
	\begin{align}\label{eq:loc dissipation intro}
		D(\dot u) =\frac12 \int_\Omega \Dd \epsi(\dot u (x)) : \epsi (\dot u(x) ) \dd x 
	\end{align}
	to be the \df{local   viscous     potential},  which    corresponds to     the nonlocal-to-local limit of $D_\rho$  defined in \eqref{eq:dissipation potential}, see Lemma \ref{lem:lim}.    The local   viscoelastic   problem \eqref{eq:pde intro}   models   the quasistatic evolution for a  linear  \name{Kelvin}-\name{Voigt} viscoelastic material in classical rheology.

	   The   energy \eqref{eq:energy} could easily be complemented with the contribution of body forces.   As this does not add difficulties, we prefer to stay with the current, force-free situation,   for the sake of  notational  brevity.  % Moreover, also nontrivial Dirichlet boundary conditions could be considered in our model with only adding notational intricacies. 

	\smallskip \smallskip 
  	\noindent\textbf{Proof strategy.}    
	Problem \eqref{eq:pde intro} is a gradient flow with respect to the metric   defined via the viscous potential as   \[ (u, \bar{u} ) \mapsto \int_\Omega \Dd \epsi(u- \bar{u}) : \epsi(u - \bar{u}). \]    We will tackle the nonlocal-to-local limit of \eqref{eq:nonlocal problem intro} as $\rho \to \delta_0$ by resorting to     an evolutionary $\Gamma$-convergence approach  \cite{mielkeGamma, sandierGammaconvergenceGradientFlows2004}.  
	  To this aim,     the gradient flow \eqref{eq:nonlocal problem intro} is reformulated in terms of an \df{energy dissipation  equality}  (also called  \df{Energy Dissipation Principle} or  \df{De Giorgi's Principle} \cite{santambrogio}),  which reads as
	\begin{align}\label{eq:edi intro}
		E_\rho\left(u_\rho(s)\right) \Big\vert_{0}^t + \int_0^t D_\rho\left(\dot u_\rho(s)\right) \dd s + \int_0^t D_\rho^* \big(-\dd_\rho E_\rho(\dot u_\rho(s))\big) \dd s = 0 \quad \forall t \in [0,T].
	\end{align}
	Here, $D^*_\rho$ denotes the convex conjugate of the   nonlocal potential of viscous forces     $D_\rho$.

	 Passing to the liminf  in equation \eqref{eq:edi intro}  and using  the inequalities
	\begin{align}
		&\liminf_{\rho \to \delta_0} E_\rho (u_\rho(t) ) \geq E(u(t)),  \label{eq:li1}\\
		&\liminf_{\rho \to \delta_0 } \int_0^t D_\rho (\dot u_\rho(s)) \dd s \geq \int_0^t D(\dot u(s))\dd s, \label{eq:li2}\\
		&\liminf_{\rho \to \delta_0 } \int_0^t D_\rho^* \big( -\dd_\rho E_\rho (u_\rho(s))\big)  \dd s \geq \int_0^t D^* \big( -\dd E(u(s))\big) \dd s \label{eq:li3}
	\end{align}
  	 as well as assuming well-preparedness of initial data (namely, $u_\rho(0) \to u(0)$ and $E_\rho(u_\rho(0)) \to E(u(0))$),     
we deduce the corresponding energy dissipation inequality for the (local) limiting problem
	\begin{align}\label{eq:edi intro local}
		E(u(t)) - E(u(0)) + \int_0^t D(\dot u(s)) \dd s + \int_0^t D^* \big(-\dd E(u(s)) \big) \dd s \leq 0.
	\end{align}
In \eqref{eq:li3} and \eqref{eq:edi intro local}, $\dd E(u)$ denotes the Fr\'echet derivative in $H^1(\Omega;\R^d)$ of the local energy $E$  \eqref{eq:loc energy intro} at $u$.  From the definition of the Fenchel-Legendre transform  and the chain rule,  we   obtain equality in \eqref{eq:edi intro local}. This is   eventually     equivalent to $u$ being a   (weak)      solution to the local   viscoelastic   problem \eqref{eq:pde intro}.  
	
	 This paper is devoted to make the above argument rigorous,   see    Theorem~\ref{thm:main} below.  In order to achieve this, we   prove      compactness and the  above   liminf-inequalities  \eqref{eq:li1}--\eqref{eq:li3}.    In fact,    inequality  \eqref{eq:li1} has already been   established     in \cite{mengeshaVariationalLimitClass2015} as part of the $\Gamma$-convergence result.   The other two liminf-inequalities    \eqref{eq:li2} and \eqref{eq:li3}     are proved in Section~\ref{sec:proof}. Crucial for proving both inequalities  is the nonlocal-to-local convergence of  $\dd_\rho D_\rho$ and $\dd_\rho E_\rho$.       Note that here also   the arguments depend on $\rho$,  which adds   some     difficulties       if    compared to   the analysis in     \cite{mengeshaVariationalLimitClass2015}.
	
	  Passing to the limit directly in the equations, namely, without reasoning to the equivalent energy-dissipation-inequality formulation, would also be possible. Still, it would require to cope with $\rho$-dependent spaces, resulting in a more involved argument.  The variational perspective has the advantage of additionally providing  convergence of the nonlocal elastic energies along the evolution,  see Corollary \ref{cor: en conv}.

		\smallskip \smallskip
  	\noindent\textbf{Literature.}    
	 Let us now put our result in the context of the available literature.  Since the introduction of peridynamics by \name{Silling} \cite{sillingReformulationElasticityTheory2000}, the dynamic peridynamic model has already attracted some attention, see \cite{ELP_surveyPD, EP_surveyDyn} and the references therein for an overview on dynamic problems and \cite{Peri_review} for a survey of the general peridynamic literature. In \cite{EmmrichWeckner_Navier}, the authors prove well-posedness of a linear peridynamic momentum balance by using the Duhamel principle. Moreover, in the nonlocal-to-local limit they recover the Navier equation of classical elasticity in a perturbative fashion, similar to that of \cite{AG_HigherOrderGradForAC, EM_CauchyBornDynamic} in the case of atomistic-to-continuum limits. This result has then been generalized in various directions, e.g., allowing for less regular solutions \cite{DZ_Ana}, or rotational variant kernels \cite{MS_solvability}. Also the nonlinear case has been studied \cite{CDMV_wellposed}.  
	
	For the static case,    in \cite{DM_bondbased} it has been shown   well-posedness of a bond-based peridynamic boundary value problem,  namely, $\mathcal{L}_\rho u = b$, where $\mathcal{L}_\rho$ is a nonlocal operator  (similar to $v\mapsto (u, v)_\rho$)  and $b$ models body forces.  Moreover,  its nonlocal-to-local limit to the classical Navier equation has been derived.  
		
	Various viscoelastic models have already been used successfully  by the engineering community to model and simulate rate-dependent phenomena involving discontinuities \cite{EnginWithComp, DLLT_visc, MO_thermoviscoelastiv, M_visco_Sandia, MW_visco, TZ_visco}. Rigorous mathematical results involving viscosity are    however    still  not available.  
	
	A nonlocal parabolic equation related to peridynamics has been investigated in \cite{Yangari}. Existence and uniqueness of solutions to the nonlocal problem are proved by a Galerkin method. However, nonlocal-to-local convergence has not been treated. In a recent work \name{Yang}, \name{Zhang}, and \name{Nie} \cite{NYZ_nonlocalDamping} considered a linearized bond-based peridynamic dynamic problem incorporating a nonlocal damping term. They establish well-posedness of the nonlocal problem by resorting to a Galerkin method, as well. Passing to the limit in the equation, they show that solutions to the nonlocal problem converge to weak solutions of the local problem, which turns out to be the  elastodynamic equation of motion. In fact,  due to their choice of the nonlocal kernel in the damping operator, the viscous effects vanish when passing to the limit and the bond-based peridynamic approach gives rise to the Laplacian in the limit (see also \cite{BCMC_noHyperEl, sillingReformulationElasticityTheory2000} for a discussion on the limitations of the bond-based approach). 
	 Although their result is similar to ours, the techniques are fundamentally different. Their result applies to a bond-based peridynamic model only   and     it does not allow to recover viscous effects in the small-horizon limit.

	A paper more closely to ours in spirit is \cite{kruzikQuasistaticElastoplasticityPeridynamics2018}. There, the authors investigate a peridynamic model related to elastoplasticity in the small strain regime. They prove well-posedness of the nonlocal problem via a time discretization scheme. Then, they pass to the limit via an evolutionary $\Gamma$-convergence approach for rate-independent systems \cite{MRS_rateIndep}. As the   viscosity     potential $D_\rho$ \eqref{eq:dissipation potential} considered here is quadratic, we are not in the setting of rate-independent systems.

		\smallskip \smallskip
	\noindent\textbf{Organization of the paper.}    
		This note is organized as follows. We introduce the model and   state  the main result on   the evolutionary $\Gamma$-convergence   in Section~\ref{sec:model and result}. In Section~\ref{sec:derivatives}, we compute the derivatives of the energy and the   viscous     potential and collect some of their properties.  These  will then be used in Section~\ref{sec:aux}, where existence of solutions to the nonlocal problem, a priori estimates, and compactness are proved. In Section~\ref{sec:conv of devs}, we collect several nonlocal-to-local convergence results and in particular establish the convergence of the derivatives of $D_\rho$ and $E_\rho$. The last   part,     Section~\ref{sec:proof}, is then devoted to the proof of our main result.

	%%%%%%%%%%%%%%%%%%%%%%%%%%%%%%%%%%%%%%%%%%%%%%%%%%%%%
	\section{The model and main result}\label{sec:model and result}
	%%%%%%%%%%%%%%%%%%%%%%%%%%%%%%%%%%%%%%%%%%%%%%%%%%%%%
	
	%%%%%%%%%%%%%%%%%%%%%%%%%%%%%%%%%%%%%%%%%%%%%%%%%%%%%
	\subsection{Notation}
	%%%%%%%%%%%%%%%%%%%%%%%%%%%%%%%%%%%%%%%%%%%%%%%%%%%%%

	We  denote the Euclidean inner product in $\R^d$ by $a \cdot b := a_i b_i$ and the contraction   of two-tensors     by  $A:B := A_{ij} B_{ij}$.
	
	Given a smooth function $u\colon \Omega \subset \R^d \to \R^d$, we abbreviate the partial derivative with respect to the $j$-th  spatial  coordinate with $\partial_j$   and write $\D u(x) \in \R^{d \times d}$   with     $(\D u(x))_{ij} = \partial_j u_i(x) $  for its Jacobian. Moreover, we denote the \df{symmetric gradient of $u$} by $\epsi(u) := (\D u + (\D u)^\top)/2$. 	  We   let  $\fint_{S^{d-1}} f(s) \dd\Hd^{d-1}(s)  := \big(\Hd^{d-1}(S^{d-1}) \big)^{-1} \int_{S^{d-1}} f(s) \dd\Hd^{d-1}(s)$, where $\Hd^{d-1}$ is the $(d-1)$-dimensional Hausdorff-measure and $S^{d-1}$ the $(d-1)$-dimensional unit sphere.     
	
	We follow the convention to abbreviate the $L^2$-inner product by $(u,v) := \int_\Omega u(x) \cdot v(x) \dd x$ for $u,v    \in L^2(\Omega;\R^d)$.   Moreover, the $L^p(\Omega;\R^N)$-norm will be denoted by $\Vert\blank \Vert_{L^p(\Omega)^N}$, where we omit  $N$  if   $N = 1$.  	Furthermore, $c$ denotes a   generic    positive constant,   only depending on data    and the dimension   $d$,        and  possibly varying from line to line. 
	
	%%%%%%%%%%%%%%%%%%%%%%%%%%%%%%%%%%%%%%%%%%%%%%%%%%%%%
	\subsection{The model}\label{subsec:model}
	%%%%%%%%%%%%%%%%%%%%%%%%%%%%%%%%%%%%%%%%%%%%%%%%%%%%%

	Let  $\Omega \subset \R^d$ ($d \geq 2$) be a bounded domain  with Lipschitz boundary    representing     the reference configuration of the elastic body. The deformation of the body is described by a \df{time-dependent displacement} $u\colon [0,T] \times \Omega \to \R^d$, where $T >0$ is a  \df{final reference time}. 
	
	To introduce  \df{nonlocal homogeneous Dirichlet boundary values} (or nonlocal volumetric constraints), see, e.g., \cite{nonlocDiff, nonlocalVectorCalc}, we fix  $\OO \subset \R^d$ open such that $\Omega \subset \subset \OO$  (compactly contained). 
	 	We then extend deformations from $\Omega$ to $\OO$ by $0$, that is, we set

		\begin{equation*}%\label{eq:V}
		V := \left\{  u \in L^2 (\OO ; \R^d) : u = 0 \textup{  on  }  \OO \setminus \Omega  \right\}.
	\end{equation*} 
	We recall  that  the nonlocal strain $\Dsc(u)$ and the nonlocal divergence $\Dfr_\rho(u)$  are given by 
	\begin{align*}
		\Dsc(u)(x',x) &= \frac{(u(x') - u(x)) \cdot (x'-x)}{| x' - x|^2},  \\
		\Dfr_\rho(u)(x) &= \pv\int_\OO \rho(x'-x) \Dsc(u)(x',x) \dd x',
	\end{align*}
	 where $\rho$ is a positive and integrable kernel,   see Assumption~\ref{assumptions} below.      Moreover, we recall   the notation         
	\begin{equation}\label{eq:defintion nonlocal seminorm}
		\vert u  \vert_{\rho}^2 = \int_{\OO} \int_{\OO} \rho(x'-x) \vert \Dsc(u) (x',x) \vert^2 \dd x' \dd x 
	\end{equation}
	 and the   definition of the     \df{nonlocal     energy space} 
	\begin{equation*}
		\S_{\rho} (\OO;\R^d) = \left\{   u \in V     :  \vert u \vert_{\rho}^2 < \infty \right\}.
	\end{equation*}
   
 Fixing the  displacement  on  $\OO\setminus \Omega$  allows us to employ a nonlocal Korn-type inequality on $ 	\S_{\rho} (\OO;\R^d)$ \cite{mengeshaNONLOCALKORNTYPECHARACTERIZATION2012}, see Proposition~\ref{prop:poincare} below, which entails that the seminorm in \eqref{eq:defintion nonlocal seminorm} is actually a norm. Therefore, $\S_{\rho}(\OO;\R^d)$	is a Hilbert space when equipped with the inner product
	\begin{equation*}
		( u,v )_{\rho} := \int_{\OO} \int_{\OO} \rho(x' - x ) \Dsc(u)(x',x) \Dsc(v)(x',x) \dd x' \dd x,
	\end{equation*}
	see \cite[Theorem~2.1]{mengeshaVariationalLimitClass2015}.

	   Recalling the definition of    the nonlocal elastic energy \eqref{eq:energy}, we have
	\begin{align}\label{eq:energy space}
		E_\rho(u) < \infty \iff \vert u \vert_{\rho} < \infty. 
	\end{align}
	Indeed, by H\"older's inequality we have that $\Dfr_\rho(u) \in L^2(\OO;\R^d)$ if  $\vert u \vert_\rho \leq c$,  and thus $E_\rho(u) \leq c$ for $u \in \S_{\rho}(\OO; \R^d)$. 	The converse also holds true: If $E_\rho(u) < \infty$, then also $\vert u \vert_{\rho} < \infty$, see \cite[Proposition~1]{mengeshaNonlocalConstrainedValue2014}. 
	 
	We extend the energy $E_\rho$ from $\S_{\rho}(\OO;\R^d)$ to all of   $V$     by setting $E_\rho(u) = \infty$ whenever $u \in   V     \setminus \S_{\rho}(\OO;\R^d)$.

	To study the passage from the nonlocal to the local model, we introduce a family of radial kernels $\rho_n :\R^d \to [0,\infty]$   suitably     converging to the Dirac $\delta_0$.  More precisely, we require $(\rho_n)_n$ to satisfy the following.

	\begin{assumption}\label{assumptions}
		Let $\rho_n : \R^d \to [0,\infty)$ ($n \geq 1$) be a family of radially symmetric kernels,   i.e.,  $\tilde{\rho}_n (\vert \xi \vert) = \rho_n(\xi)$ for some $\tilde{\rho}_n\colon [0,\infty) \to [0,\infty)$.   
		We assume that $(\rho_n)_n$ satisfies the following
		\begin{alignat}{2}
			\mathrm{(i)}& \quad && r\mapsto r^{-2} \tilde{\rho}_n (r)  \text{ is nonincreasing in } r,  \nonumber \\
			\mathrm{(ii)}& &&\int_{\R^d} \rho_n ( \xi ) \dd  \xi   = d,     \label{eq:kernel normalization}\\
			\mathrm{(iii)}& &&\lim_{n \to \infty} \int_{  \{ \vert  \xi  \vert \geq r_0 \} } \rho_n ( \xi ) \dd  \xi  = 0 \qquad \forall r_0 >0. \nonumber
		\end{alignat}
	\end{assumption} 
 These assumptions have already been considered in \cite{mengeshaVariationalLimitClass2015}   for the proof of the $\Gamma$-convergence  of  $(E_{\rho_n})_n$.    
  In the following, in order to simply notation,  we set \[  E_n := E_{\rho_n}, \qquad  D_n := D_{\rho_n}, \qquad  \Sn(\OO;\R^d) := \S_{\rho_n}(\OO;\R^d) , \qquad \Dfr_n := \Dfr_{\rho_n}, \qquad    \vert \cdot \vert_n := \vert \cdot \vert_{\rho_n}.     \]  
 As  $n \to \infty$,   one expects the   nonlocal     energy spaces $\Sn(\OO;\R^d)$ to approach    the set $\big\{ u \in H^1(\OO; \R^d) : u = 0 \textup{  on  } \OO\setminus \Omega \big\}$, which can be identified with       
\begin{align*}
	H^1_0(\Omega;\R^d) := \left\{ H^1(\Omega;\R^d) : u\vert_{\partial\Omega}  = 0    \right\},
\end{align*}
	where $u\vert_{\partial\Omega}$ has to be interpreted in the   usual     trace sense.  %The space $H^1_0(\Omega;\R^d)$ is again Hilbert, as a Korn-inequality holds, see Lemma~\ref{lem:korn type} below.  

	%%%%%%%%%%%%%%%%%%%%%%%%%%%%%%%%%%%%%%%%%%%%%%%%%%%%%
	\subsection{Main result}
	%%%%%%%%%%%%%%%%%%%%%%%%%%%%%%%%%%%%%%%%%%%%%%%%%%%%%
	Given some \df{initial nonlocal displacement} $u_{n}^0 \in  \Sn(\OO;\R^d)    $, we   consider    the     solution     $u_n$    to the \df{nonlocal   viscoelastic    problem} given by
	\begin{subequations}	  \label{eq:nonlocal problem}
		\begin{empheq}[left=\empheqlbrace]{alignat = 2}
		 	\dd_n  D_n (\dot u_n) + \dd_n E_n (u_n) &= 0  &\quad  &\text{in }   \Sn(\OO;\R^d)^*  \quad   \textup{ for  almost all } t \in [0,T],         \label{eq:nonlocal problem 1} \\
			u_n(0,\cdot) &= u_{n}^0  &\quad &\text{in } \OO. 
		\end{empheq}  
	\end{subequations}

	Here, $\dd_n F(x)$ denotes the   Fr\'echet differential  of $F$ evaluated at $x$  in the space $\Sn(\OO;\R^d)$   and $(\Sn(\OO;\R^d))^*$ denotes the dual space of $\Sn(\OO;\R^d)$.     We say that $u_n \in   H^1    (0,T;\Sn(\OO;\R^d))$ is a solution to \eqref{eq:nonlocal problem} if it satisfies the initial condition and
	\begin{align*}
	 	\dd_n D_n (\dot u_n)(v) + \dd_n E_n( u_n) (v) = 0
	\end{align*}
	 for all $v \in \Sn (\OO;\R^d)$ and    almost every $t \in [0,T]$,  where     $\dd_n F (x)(v) = \lim_{\eta \to 0} \frac{F(x+\eta v) - F(x) }{\eta}$ is the G\^{a}teaux derivative of $F$ at $x \in \Sn(\OO;\R^d)$ in direction $v \in \Sn (\OO;\R^d)$.

	The \df{local   viscoelastic   problem} is to find $u \in    H^1    (0,T;H^1_0(\Omega;\R^d))$ solving,   in the weak sense, the system    
	\begin{subequations}\label{eq:pde}
		\begin{empheq}[left=\empheqlbrace]{alignat =2 }
	 	\dd D(\dot u(t)) + \dd E(u(t))  &= 0 &\quad &  \textup{in  } (H^1_0(\Omega;\R^d))^*  \quad   \textup{  for  almost every } t \in [0,T],       \\
			 	u(0,\cdot) &= u^0 &\quad &\textup{in  } \Omega.
		\end{empheq}  
	\end{subequations}
	Here,   $u^0 \in H^1_0(\Omega;\R^d)$     is a given \df{initial local displacement}.   Moreover, $\dd F (x)$ denotes the Fr\'echet differential of $F$ evaluated at $x$ in the space $H^1_0 (\Omega;\R^d)$. Owing to the explicit form of $\dd E$ and $\dd D$, see \eqref{eq:gateaux dE} and \eqref{eq:gateaux dD} below, a weak solution $u \in   H^1   (0,T;H^1_0(\Omega;\R^d))$ to \eqref{eq:pde} equivalently satisfies
	\begin{align*}
		\int_\Omega  \Dd \epsi(\dot u (x)) : \epsi (v(x) ) + \C \epsi(u(x)) : \epsi (v(x) ) \dd x = 0  
	\end{align*}
	for all $v \in H^1_0(\Omega;\R^d)$ and  almost  every $t \in [0,T]$.     Here,   recall that $\C$ is the   elasticity     tensor and $\Dd$ is the   viscosity     tensor defined in \eqref{eq:tensors intro} and \eqref{eq:tensor intro D}, respectively.

 Our main result   reads     as follows.  
	
	\begin{theorem}[Evolutionary $\Gamma$-convergence]\label{thm:main}
		Under Assumptions~\textup{\ref{assumptions}}, let a sequence of initial conditions  $u_n^0 \in \Sn(\OO;\R^d)$ and  $u^0 \in H^1_0(\Omega;\R^d)$     be given such that $u_{n}^0 \to u^0$ strongly in $L^2(\OO;\R^d)$ and $E_n(u_n^0) \to E (u^0)$ as $n \to \infty$.
		Then, for all $n \in \N$, the nonlocal   viscoelastic     problem \eqref{eq:nonlocal problem} admits a unique solution $u_n \in   H^1    (0,T; \Sn (\OO; \R^d))$ and $u_n   \to u$  in  $C([0,T]; L^2(\OO;\R^d))$,  where $u \in    H^1    (0,T;H^1_0(\Omega;\R^d))$ is   the unique weak     solution to   the local viscoelastic problem    \eqref{eq:pde}.
	\end{theorem}

  	The convergence of $u_n \to u$ and $u_n^0 \to u^0$ in the theorem above, has to be interpreted as convergence of  $u_n$ and $u_n^0$  to the extensions of $u$ and $u^0$ to $\OO$ by $0$, i.e., as convergence on $\OO$.    From now on, we will understand all limits  of elements of $\Sn(\OO,\R^d)$  in this sense without changing the notation,   and identify $u \in H^1_0(\Omega;\R^d)$ with its trivial extension to $\OO$, whenever convenient.

	\begin{remark}[Regularity]\label{rem:time}
{\rm (i)} 	As we are in the force-free setting, the solution $u$ to \eqref{eq:pde} actually satisfies $u \in C^\infty(0,T; H^1_0(\Omega;\R^d))$. Indeed, as $u \in H^1(0,T;H^1_0(\Omega;\R^d))$, also $\dd E(u) \in H^1(0,T;(H^1_0(\Omega;\R^d))^*)$ by \eqref{eq:gateaux dE} below. Comparison in equation~\eqref{eq:pde} yields that also $\dd D (\dot u ) \in H^1(0,T;(H^1_0(\Omega;\R^d))^*)$, and hence by \eqref{eq:gateaux dD} and Korn's inequality   (see  Lemma~\textup{\ref{lem:korn type}} below) we deduce $\dot u \in H^1(0,T;H^1_0(\Omega;\R^d))$. An induction argument yields the claim. 
		
{\rm (ii)}  Analogously to \textup{(i)}, if $u_n$ is a solution to the nonlocal viscoelastic problem \eqref{eq:nonlocal problem}, then  $u_n \in C^\infty(0,T; \Sn(\OO;\R^d))$ (replacing \eqref{eq:gateaux dE}--\eqref{eq:gateaux dD} by \eqref{eq:derivative of energy}--\eqref{eq:derivative of dissipation} in the argument above).
				
{\rm (iii)} 	Note that in the   case with body forces  $f \in L^2(0,T;L^2(\Omega;\R^d))$, the solution to the local viscoelastic problem \eqref{eq:pde} satisfies  (only)  $u \in H^1(0,T;H^1_0(\Omega;\R^d))$ and analogously a solution to the nonlocal viscoelastic problem would be in $H^1(0,T;\Sn(\OO;\R^d))$. 
	\end{remark}

	%%%%%%%%%%%%%%%%%%%%%%%%%%%%%%%%%%%%%%%%%%%%%%%%%%%%%
	\section{Differentials of the potentials}\label{sec:derivatives}
	%%%%%%%%%%%%%%%%%%%%%%%%%%%%%%%%%%%%%%%%%%%%%%%%%%%%%
	
	\subsection{Nonlocal potentials}
	In this section, we calculate the directional derivatives of the   nonlocal elastic     energy \eqref{eq:energy} and the   nonlocal viscous    potential \eqref{eq:dissipation potential}. 	
	  Based on this, we also get that     $D_n$ and $E_n$ are Fr\'echet differentiable, as the G\^{a}teaux derivatives are linear and continuous. 
	
	\begin{lemma}[G\^{a}teaux derivatives]\label{lem:gateaux}
		Let $u,    \, v     \in \Sn(\OO;\R^d)$. %and $v \in C^\infty_{\textup{c}}(\OO;\R^d)$. 
		Then, 
		\begin{align}
			\dd_n E_{n}(u)(v)   &= \beta  \int_\OO ( \Dfr_n(u)(x)) ( \Dfr_n(v)(x))  \dd x 
			+ \alpha  \int_\OO \int_\OO \rho_n(x'-x) \notag \\ 
			&\times\left(  \Dsc(u)(x',x) - \frac1d \Dfr_n(u)(x)  \right) \left(  \Dsc(v)(x',x) - \frac1d \Dfr_n(v)(x)  \right) \dd x' \dd x, \label{eq:derivative of energy} \\
		 \dd_n D_n (u)(v)&=  \int_\OO \int_\OO \rho_n(x' - x) \Dsc(u )(x',x) \Dsc( v) (x', x) \dd x' \dd x. \label{eq:derivative of dissipation} 
		\end{align}
	 In particular, the G\^{a}teaux derivatives are linear.  
	\end{lemma}
	\begin{proof}
		We split the energy $E_n$ into the   dilatational     part and the   deviatoric     part as   $E_n  =  E_\beta  + E_\alpha $ with  
		\begin{align*}
			E_\beta (u) &= 	 \frac{\beta}{2} \int_\OO (\Dfr_n(u)(x))^2 \dd x, \\
			E_\alpha(u) &= \frac{\alpha}{2} \int_\OO \int_\OO \rho_n(x'-x) \left(\Dsc(u)(x',x) - \frac1d \, \Dfr_n(u)(x)\right)^2 \dd x' \dd x. 
		\end{align*}
		Recall the   definition of the  nonlocal strain from \eqref{eq:nonlocal strain}.  As $u \mapsto \Dsc(u)$ is linear,  its G\^{a}teaux  derivative in direction $v \in \Sn(\OO;\R^d)$ is 
	 	\begin{align}\label{eq:deriv nonlocal strain}
			\frac{\dd}{\dd h} \Dsc (u+hv)(x',x) \big\vert_{h=0} =   \Dsc(v)(x',x).
		\end{align} 
		By \eqref{eq:deriv nonlocal strain}, for the nonlocal divergence \eqref{eq:nonlocal divergence} we get
	 	\begin{align*}
			\frac{\dd}{\dd h} \Dfr_n(u+hv) (x) \vert_{h = 0}  = \pv\int_\OO \rho_n(x'-x) \Dsc(v)(x',x) \dd x' =  \Dfr_n(v)(x), 
		\end{align*}  and consequently
	 	\begin{align*}
			\frac{\dd}{\dd h} E_\beta (u + h v) \vert_{h = 0} &=   \frac{\beta}{2}     \int_\OO	\frac{\dd}{\dd h} (\Dfr_n(u+hv)(x))^2 \big\vert_{h=0} \dd x  =  \beta    \int_\OO ( \Dfr_n(u)(x)) ( \Dfr_n(v)(x))  \dd x. 
		\end{align*}
		For the   deviatoric     energy  term,   we compute
		\begin{align*}
			&\frac{\dd}{\dd h}  E_\alpha(u + hv) \big\vert_{h=0}\\
			&= \frac{\dd}{\dd h} \frac{\alpha}{2}   \int_\OO     \int_\OO \rho_n (x'-x) \left( \Dsc(u + hv)(x',x) - \frac1d \Dfr_n(u+hv) (x) \right)^2 \dd x' \dd x \big\vert_{h=0} \\
			&=  \alpha   \int_\OO    \int_\OO \rho_n(x'-x) \left(  \Dsc(u)(x',x) - \frac1d \Dfr_n(u)(x)  \right) \left( \Dsc(v)(x',x) - \frac1d \Dfr_n(v)(x)  \right) \dd x' \dd x,
		\end{align*}
		 giving  \eqref{eq:derivative of energy}.

	 	Finally,   we calculate the  directional derivative of   $D_n$ at $ u$ in direction $v$. By \eqref{eq:deriv nonlocal strain} we readily get
		\begin{align*}%\label{eq:var of nonloc diss}
		\frac{\dd}{\dd h} D_n( u + h  v) \big \vert_{h=0} =  \int_\OO \int_\OO \rho_n(x' - x) \Dsc(u )(x',x) \Dsc( v) (x', x) \dd x' \dd x, 
		\end{align*}
		which concludes the proof.
	\end{proof}

	Recall the definition of the nonlocal seminorm~\eqref{eq:defintion nonlocal seminorm}.
	\begin{corollary}\label{cor:frechet}
		For $u, v \in \Sn(\OO;\R^d)$, we have
		\begin{align*}
			&\vert \dd_n E_n(u)(v)\vert \leq c \vert u \vert_{n}  \vert v \vert_{n}, \\
			&\vert \dd_n D_n (u)(v) \vert \leq c \vert u \vert_{n}  \vert v \vert_{n}. 
		\end{align*}
		 In particular, $E_n$ and $D_n$ are Fr\'echet differentiable in $\Sn(\OO;\R^d)$.  
	\end{corollary}
	\begin{proof}
		 Linearity is clear from Lemma~\ref{lem:gateaux}, so we focus on boundedness.	  Writing $\rho_n = \sqrt{\rho_n} \sqrt{\rho_n}$,     by H\"older's inequality we have
		\begin{align*}
			\vert \Dfr_n (u)(x) \vert^2 \leq \int_\OO \rho_n(x' - x ) \dd x' \, \int_\OO \rho_n(x' - x) \left(\Dsc(u)(x',x) \right)^2 \dd x',
		\end{align*}
		and consequently,  	 \eqref{eq:kernel normalization} yields 
		\begin{align*}
			\Vert \Dfr_n(u) \Vert_{L^2(\OO)} \leq    \Vert \rho_n \Vert_{L^1(\R^d)}^{ 1/2 }   \vert u \vert_{n}   = d^{1/2}  \vert u \vert_{n}.    
		\end{align*}
  In a similar fashion,  we also get    
\begin{align}\label{eq:bound on Dfr2}
			\int_\OO ( \Dfr_n(u)(x)) ( \Dfr_n(v)(x))  \dd x  \leq     \Vert \rho_n \Vert_{L^1(\R^d)}^{ 1/2 }  \vert u \vert_{n}  \vert v \vert_{n} =   d    \vert u \vert_{n}  \vert v \vert_{n}. 
		\end{align}
		Expanding the product in 
		\begin{align*}
			\int_\OO \int_\OO \rho_n(x'-x) \left(  \Dsc(u)(x',x) - \frac1d \Dfr_n(u)(x)  \right)\left(  \Dsc(v)(x',x) - \frac1d \Dfr_n(v)(x)  \right) \dd x' \dd x, 
		\end{align*}
		and noting that the mixed terms are of the form   $\int_\OO \Dfr_n(u)(x) \Dfr_n(v)(x) \dd x$,     one realizes that,   due to \eqref{eq:bound on Dfr2},     we just need to bound the term  
		\begin{align*}
				\int_\OO \int_\OO \rho_n(x'-x) \Dsc(u)(x',x) \Dsc(v)(x',x) \dd x' \dd x. 
		\end{align*}
		Applying H\"older's inequality, we have
		\begin{align*}
				&\int_\OO \rho_n(x'-x) \Dsc(u)(x',x) \Dsc(v)(x',x) \dd x'\\ 
				&\leq  \left(	\int_\OO \rho_n(x'-x)\left( \Dsc(u)(x',x)\right)^2 \dd x' \right)^{1/2}\left( \int_\OO \rho_n(x'-x)\left( \Dsc(v)(x',x)\right)^2 \dd x' \right)^{1/2}
		\end{align*}
		and therefore also
		\begin{equation*}
				\int_\OO \int_\OO \rho_n(x'-x) \Dsc(u)(x',x) \Dsc(v)(x',x) \dd x' \dd x \leq c \vert u\vert_{n}\vert v \vert_{n}.\qedhere
		\end{equation*}
	\end{proof}

	\subsection{Local potentials}
	
	Let us recall the definition of the tensors \eqref{eq:tensors intro} and \eqref{eq:tensor intro D}, as well as the   local elastic    energy \eqref{eq:loc energy intro} and the    local     viscous     potential \eqref{eq:loc dissipation intro}. As the tensors $\C$ and $\Dd$ are major-symmetric, i.e., $\C_{ijkl} = \C_{klij}$ and $\Dd_{ijkl} = \Dd_{klij}$, we can  compute the G\^{a}teaux derivatives of $E$ and $D$ in $H^1(\Omega;\R^d)$ as
		\begin{align}
				 \dd E(u)(v)  =  \int_\Omega \C \epsi(u(x)): \epsi(v(x)) \dd x, \label{eq:gateaux dE}\\
				 \dd D(u)(v) =  \int_\Omega \Dd \epsi(u(x)): \epsi(v(x)) \dd x,  \label{eq:gateaux dD}
		\end{align}
		for    all       $u, v \in H^1(\Omega;\R^d)$.

	%%%%%%%%%%%%%%%%%%%%%%%%%%%%%%%%%%%%%%%%%%%%%%%%%%%%%
	\section{Nonlocal well-posedness and energy-dissipation-equality}\label{sec:aux}
	%%%%%%%%%%%%%%%%%%%%%%%%%%%%%%%%%%%%%%%%%%%%%%%%%%%%%
	
	In this  section, we prove    the    well-posedness of the nonlocal   viscoelastic    problem~\eqref{eq:nonlocal problem}, derive the energy-dissipation equality, and prove compactness. To this end, we start with a short preliminary section,  recalling   some     related  results from the literature.

	%%%%%%%%%%%%%%%%%%%%%%%%%%%%%%%%%%%%%%%%%%%%%%%%%%%%%
	\subsection{Preliminary results}
	%%%%%%%%%%%%%%%%%%%%%%%%%%%%%%%%%%%%%%%%%%%%%%%%%%%%%
	
	 To start with, let us recall the following useful result.
	\begin{lemma}[{\cite[Lemma~2.1]{mengeshaNONLOCALKORNTYPECHARACTERIZATION2012}}]\label{lem:bound Sn seminorm by H1}
		For any $u \in H^1(\OO;\R^d)$, we have
		\begin{align*}
			\vert u \vert_n \leq   c  \Vert \epsi(u) \Vert_{L^2(\OO)^{\dxd}},
		\end{align*}    
		  where $c$ only depends on $\OO$ and $d$.    
	\end{lemma}
	Next, we consider a nonlocal Poincar\'e-Korn inequality. To this end, note that $\vert u \vert_{n} = 0$  if and only if   $u$ is an infinitesimal rigid displacement, i.e.,
	\begin{equation*}
		u\in	\Rigid := \left\{ x \mapsto Ax + t :  A \in \R^{d\times d}, \ A^\top = -A, \ t \in \R^d   \right\},
	\end{equation*}
	see \cite[Lemma~1]{mengeshaNonlocalConstrainedValue2014}.   In particular,  $\Sn(\OO;\R^d)  \cap \Rigid = \{0 \}$, which is necessary for the following.      

	\begin{proposition}[Nonlocal Poincar\'{e}-Korn inequality, {\cite[Proposition~2.7]{mengeshaVariationalLimitClass2015}}]\label{prop:poincare}
		There exists a positive constant   $c = c(\OO, \rho_n)$     such that 
		\begin{equation}\label{eq:korn}
			\Vert u \Vert_{L^2(\OO)^d}^2 \leq c \vert u \vert_{n}^2 \qquad \forall u \in   \Sn(\OO;\R^d).   
		\end{equation}
	\end{proposition} 
  Note that in principle  the constant may depend on $n$. Yet, the following compactness result and a standard compactness-contradiction argument   show  that  the constant in  \eqref{eq:korn} can    in fact    be chosen independently of $n$.

	 \begin{proposition}[Compactness, {\cite[Proposition~4.2]{mengeshaVariationalLimitClass2015}}]\label{prop:compactness static}
	 	Assume that $\rho_n$    satisfy     Assumption~\textup{\ref{assumptions}}. If $u_n$ is a         bounded  sequence  in   $V$     satisfying the uniform bound 
	 	\begin{align*}
	 		\sup_{n \in \N} \vert u_n \vert_{n} \leq c,
	 	\end{align*}
	 	there exists    $u \in H^1_0(\Omega;\R^d)$     and a (not relabeled) subsequence such that $u_n \to u$ strongly in $L^2(\OO;\R^d)$. 
	 \end{proposition}    
 
  To be precise, the fact that $u$ has zero trace is not inferred from \cite[Proposition~4.2]{mengeshaVariationalLimitClass2015}, but  comes from the almost everywhere pointwise convergence of   $u_n \in V$,   which is $0$ on $\OO \setminus \Omega$.  
	   
	 In particular,   the estimate \eqref{eq:korn}    entails that $\vert \cdot \vert_n$ is a norm on $\Sn(\OO;\R^d)$ and thus $( \cdot , \cdot )_n$ constitutes an inner product. 
	Note that relation \eqref{eq:derivative of dissipation} implies that $\dd_n D_n$ is the Riesz isomorphism, as $\dd_n D_n(u)(v) = (u,v)_n$.

		We recall the classical Korn inequality,  see, e.g., \cite[Proposition~1]{FrieseckeJamesMueller}.    
	
\begin{lemma}[Korn inequality]\label{lem:korn type}
	There is a constant $c > 0$ with 
	\begin{align*}
		\Vert u \Vert_{H^1(\Omega)^d} \leq c \Vert \epsi(u) \Vert_{L^2(\Omega)^{\dxd}} \qquad \forall u \in H^1_0(\Omega;\R^d).
	\end{align*}
\end{lemma}

 Next, we recall  the following $\Gamma$-convergence result, which is crucial for our analysis.

\begin{theorem}[$\Gamma$-convergence,  {\cite[Theorem~1.2]{mengeshaVariationalLimitClass2015}}]\label{thm:gamma convergence}
	The sequence $E_n$ $\Gamma$-converges in the strong $L^2$-topology to $E$,  defined in \eqref{eq:gamma limit of energy},  which can be equivalently rewritten as     
	\begin{align}\label{eq:energy in gamma conv thm}
		E (u) =&\frac12 \int_\Omega \C \epsi(u(x)): \epsi(u(x)) \dd  x \notag \\
		=& \frac{\beta}{2} \int_\Omega \div^2 u(x) \dd x  + \frac{\alpha}2   d    \int_\Omega \fint_{S^{d-1}} \left(s \cdot \D u (x) s - \frac1d \div u(x) \right)^2 \dd \Hd^{d-1}(s) \dd x \notag  \\
		=& \frac{\mu}{2} \int_\Omega \vert \epsi(u(x)) \vert^2 \dd x	+ \frac{\lambda}{2} \int_\Omega (\div u(x))^2 \dd x,
	\end{align}
	with   $\mu = \frac{2\alpha}{d+2}$     and   $\lambda =\beta - \frac{2\alpha}{d (d+2)}$.     
\end{theorem}

A proof of the claim that the different representations   \eqref{eq:loc energy intro} and \eqref{eq:energy in gamma conv thm}     of the energy $E$ are indeed   equivalent     can be found in the Appendix~\ref{sec:appendix calculations}, Lemmas~\ref{lem:energy rewritten} and \ref{lem:energy as tensor}.

 We close this section with a representation of the viscous potential $D$ which  follows from Lemma~\ref{lem:prod  sDus} and a computation analogous to the one of Lemma~\ref{lem:energy as tensor}  (by just replacing $\lambda$ and $\mu$ there by  $\frac{1}{d+2}$ and $\frac{2}{d+2}$, respectively). 
 
   \begin{lemma}\label{lem:vis potential forms}
 	The   local     viscous potential, defined in \eqref{eq:loc dissipation intro} can be equivalently written as
 	\begin{align*}
 		D(w) &= \frac12 \int_\Omega \Dd \epsi (w(x)) : \epsi (w(x)) \dd x \\
 		&= \frac{1}{d+2} \int_\Omega \vert \epsi (w(x)) \vert^2 \dd x + \frac{1}{2(d+2) } \int_\Omega (\div w(x))^2 \dd x \\
 		&= \frac{d}{2} \int_\Omega \fint_{S^{d-1}} (s \cdot \D w (x) s)^2 \dd \Hd^{d-1}(s) \dd x
 	\end{align*}
 	for  all  $w \in H^1(\Omega;\R^d)$.
 \end{lemma}

	%%%%%%%%%%%%%%%%%%%%%%%%%%%%%%%%%%%%%%%%%%%%%%%%%%%%%
		\subsection{Solutions to the nonlocal   viscoelastic   problem}\label{subsec:proof solution to nonlocal}
	%%%%%%%%%%%%%%%%%%%%%%%%%%%%%%%%%%%%%%%%%%%%%%%%%%%%%

	Here, we prove   the     well-posedness of the nonlocal   viscoelastic   problem~\eqref{eq:nonlocal problem}. As $\dd_n D_n$ and $\dd_n E_n$ are linear, the nonlocal   viscoelastic   problem is nothing but a linear Cauchy problem on the Hilbert space $\Sn (\OO;\R^d)$. Therefore, well-posedness is classical. %For the sake of completeness, we nevertheless present a proof. 
	
	\begin{proposition}[Well-posedness of nonlocal   viscoelastic    problem]\label{prop:wellposedness of nonlocal problem}
		The nonlocal   viscoelastic    problem~\eqref{eq:nonlocal problem} has a unique solution $u_n \in   H^1    (0,T;\Sn(\OO; \R^d))$.
	\end{proposition}	
	\begin{proof}
		As   $\dd_n D_n(\dot u_n)(v_n) = (\dot u_n, v_n)_n$     by Lemma~\ref{lem:gateaux},  $\dd_n D_n$ is nothing but the Riesz isomorphism. In particular, \eqref{eq:nonlocal problem} is equivalent to 
		\begin{align*}
			\frac{\dd}{\dd t} u_n (t) + ((\dd_n D_n)^{-1} \circ \dd_n E_n)(u_n) = 0.
		\end{align*}
		 Moreover,  as   $\dd_n D_n$ is  the Riesz isomorphism, Corollary~\ref{cor:frechet} implies that $(\dd_n D_n)^{-1} \circ \dd_n E_n$ is a Lipschitz map. Therefore, the classical Picard Theorem for Banach spaces applies, see, e.g., \cite[Theorem~6]{Lobanov_ODE}. 
\end{proof}

	%%%%%%%%%%%%%%%%%%%%%%%%%%%%%%%%%%%%%%%%%%%%%%%%%%%%%
	\subsection{Energy dissipation   equality     and a priori estimates}
	%%%%%%%%%%%%%%%%%%%%%%%%%%%%%%%%%%%%%%%%%%%%%%%%%%%%%
	
	In this section, we derive the energy dissipation equality and, as a consequence,  the   a priori estimates. As a crucial ingredient, we first check that the chain rule holds in the nonlocal   setting.

	\begin{lemma}[Chain rule]\label{lem:nonloc chain rule}
		Let $u_n$ be the solution to the nonlocal   viscoelastic    problem \eqref{eq:nonlocal problem}.
	   Then, the chain rule holds, i.e., we have, for almost all $t \in (0,T)$,    
		\begin{align*}
			\frac{\dd}{\dd t} E_n (u_n(t, \cdot)) =    \dd_n E_n (u_n(t)) (\dot u_n (t) ).    
		\end{align*}
	\end{lemma}
	\begin{proof}
		By Proposition~\ref{prop:wellposedness of nonlocal problem}, we have that $\dot u_n \in L^2(0,T;\Sn(\OO;\R^d))$. Corollary~\ref{cor:frechet} implies that $\dd_n D_n (\dot u_n) \in L^2(0,T;(\Sn(\OO;\R^d))^*)$  and  hence by comparison in  equation~\eqref{eq:nonlocal problem}, also $\dd_n E_n (u_n) \in L^2(0,T;(\Sn(\OO;\R^d))^*)$. The assertion now follows from \cite[Proposition~2.2]{ulisse},   after noting that $E_n$ is quadratic and therefore also proper, convex, and lower semicontinuous.  
	\end{proof}

  	 	Let $u_n$ be the solution to the nonlocal   viscoelastic   problem~\eqref{eq:nonlocal problem}, see Proposition~\ref{prop:wellposedness of nonlocal problem}. Testing the nonlocal   viscoelastic   problem \eqref{eq:nonlocal problem} against $\dot{u}_n$, we obtain    
	\begin{align*}
		  0 =  \dd_n E_n(u_n) (\dot{u}_n )+  \dd_n D_n(\dot{u}_n)(\dot{u}_n).    
	\end{align*}
	  As equality in   the     Fenchel-Young inequality holds, using the chain rule Lemma~\ref{lem:nonloc chain rule}, integration in time yields, for all $t \in [0,T]$,  
	\begin{align}\label{eq:edi}
		E_n(u_n(s, \cdot)) \Big\vert_{0}^t + \int_0^t D_n(\dot u_n(s, \cdot)) \dd s + \int_0^t D_n^* \big( -\dd_n E_n(u_n(s, \cdot)) \big) \dd s = 0,
	\end{align}
	where   $D_n^*$     denotes the Fenchel-Legendre transform of $D_n$,   namely, 
	\begin{align}\label{eq:defin D*}
		D_n^*(v^*) = \sup \left\{ \langle v^* , w \rangle_n - D_n(w) : w \in \Sn (\OO;\R^d) \right\}.
	\end{align}
 Here,  $\langle \cdot , \cdot \rangle_n$ denotes the dual pairing  between $\Sn(\OO;\R^d)^*$ and        $\Sn (\OO;\R^d)$.  From \eqref{eq:edi} and the assumption that $\sup_n E_n(u_n(0))< \infty$,    see  the  statement of Theorem~\ref{thm:main},    we immediately get the following lemma.
	\begin{lemma}[A priori estimates]\label{lem:a priori estimates}
		Let $u_n$   be the     solution to the nonlocal   viscoelastic    problem~\textup{\eqref{eq:nonlocal problem}}. Then, for  all $t \in [0,T]$ we have 
		\begin{align}
			&E_n(u_n(t, \cdot)) \leq c, \label{eq:unif E}\\
			&\int_0^t D_n(\dot{u}_n(s, \cdot) ) \dd s \leq  c, \label{eq:unif D} \\
			&\int_0^t D^*_n \big(- \dd_n E_n (u_n(s, \cdot ) ) \big) \dd s \leq c, \label{eq:unif D stern}
		\end{align}
		where $c$ only depends on $\sup_n E_n(u_n(0,\cdot))$.
	\end{lemma}

 Here, we use  that $D^*_n  ( v_n^* )\geq 0$ for all $v_n^* \in (\Sn(\OO;\R^d))^*$   as $D_n(0) = 0$.      Furthermore,  by   the     definition of $D_n$, \eqref{eq:unif D} is nothing else but
	\begin{align}\label{eq:unif bound dot u}
		\int_0^t \vert \dot u_n (s, \cdot ) \vert_{n}^2 \dd s \leq c.
	\end{align}

	%%%%%%%%%%%%%%%%%%%%%%%%%%%%%%%%%%%%%%%%%%%%%%%%%%%%%
	\subsection{Compactness}
	%%%%%%%%%%%%%%%%%%%%%%%%%%%%%%%%%%%%%%%%%%%%%%%%%%%%%
	
	    Recall the static compactness result Proposition~\ref{prop:compactness static}, which  allows us to prove an Aubin-Lions-type result.
	
	\begin{proposition}[Aubin-Lions-type compactness]\label{prop:compactness}
		Let $(u_n)_n$ be a sequence of solutions to the nonlocal   viscoelastic   problem \eqref{eq:nonlocal problem}.
		Then, there exists $u \in C([0,T]; L^2( \OO; \R^d) )$  and a (non relabeled) subsequence $(u_n)_n$ such that 
		\begin{equation*}
			u_n \to u \quad \text{in  } C([0,T];L^2(\OO;\R^d)).
		\end{equation*}
	\end{proposition}
	\begin{proof}
		We use the Arzel\`{a}-Ascoli Theorem. 
		First, note that    \eqref{eq:energy space} and     \eqref{eq:unif E}  imply  that the set $\{ u_n(t) : n \in \N \}$   is   bounded     in $\Sn(\OO;\R^d)$. Thus,   by Proposition~\ref{prop:compactness static}  it is relatively compact in    $L^2( \OO; \R^d)$. % embedding is continuous and continuous composed with compact is compact
		Equicontinuity follows from \eqref{eq:unif D}. Indeed, from   \eqref{eq:unif bound dot u}     and the nonlocal Poincar\'{e} inequality, Proposition~\ref{prop:poincare}, we have that $\sup_n \Vert \dot{u}_n \Vert_{L^2(0,T;L^2(\OO;\R^d))} \leq c$. Thus, 
		\begin{align*}
			\Vert u_n(t) - u_n(s) \Vert_{L^2(\OO)^d} \leq \int_s^t \Vert \dot{u}_n(\tau) \Vert_{L^2(\OO)^d} \dd \tau \leq \sqrt{t-s} \Vert \dot{u}_n \Vert_{L^2(0,T;L^2(\OO;\R^d))} \leq c \sqrt{t-s}.
		\end{align*}
		Therefore, the Arzel\`{a}-Ascoli Theorem is applicable and yields the claim.
	\end{proof}
   Note that  $u(t) \in H^1_0(\Omega,\R^d)$ for all $t \in [0,T]$ from Proposition~\ref{prop:compactness static}.  
	
	%%%%%%%%%%%%%%%%%%%%%%%%%%%%%%%%%%%%%%%%%%%%%%%%%%%%%
	\section{Nonlocal-to-local convergence of derivatives}\label{sec:conv of devs}
	%%%%%%%%%%%%%%%%%%%%%%%%%%%%%%%%%%%%%%%%%%%%%%%%%%%%%

	\subsection{General nonlocal-to-local convergence results}\label{subsec:nonloc conv}
	
	In this section, we collect several results,   which will be used in the following.

	\begin{lem}[Divergence]\label{lem:convergence of strain}
		For   $u \in H^1_0(\Omega; \R^d)$     we have 
		\begin{align*}
			\pv\int_\OO \rho_n (x'-\cdot) \Dsc(u)(x', \cdot) \dd x' \to   d     \fint_{S^{d-1}} s\cdot \D u s \dd\Hd^{d-1}(s) =   \div u    
		\end{align*}
		strongly in $L^2 (\OO;\R^d) $.
	\end{lem}
   A proof  of the above lemma  can be found in \cite[Lemma 3.1]{mengeshaVariationalLimitClass2015},    where 	  the second identity  follows from Lemma~\ref{lem:div}.

	 A similar result holds true if  the function $u$ is replaced by a sequence of functions $(u_n)_n$,   up to resorting to     weak convergence. 
	 \begin{lem}[{\cite[Lemma~3.6]{mengeshaVariationalLimitClass2015}}]\label{lem:convergence of strain for u delta}
	 	For   $u \in H^1_0(\Omega; \R^d)$     and a sequence $u_n \in \Sn(\OO;\R^d)$ such that $u_n \to u$ in $L^2(\OO;\R^d)$ and $\sup_n \vert u_n \vert_n \leq c$ we have 
	 	\begin{align*}
	 		\pv\int_\OO \rho_n (x'-\cdot ) \Dsc(u_n)(x', \cdot) \dd x' \wto   d     \fint_{S^{d-1}} s\cdot \D u s \dd\Hd^{d-1}(s) 
	 	\end{align*}
	 	weakly in $L^2(\OO;\R^d)$. 
	 \end{lem}
	 
	   The proof in \cite{mengeshaVariationalLimitClass2015}     is based on a nonlocal integration by parts formula.  Note that Lemmas~\ref{lem:convergence of strain} and \ref{lem:convergence of strain for u delta} immediately imply the following result.
	 \begin{corollary}\label{cor:convergece of prod of div}
	 	Let $u, v \in   H^1_0(\Omega;\R^d)    $ and $u_n \in \Sn(\OO;\R^d)$ with $u_n \to u$  in $L^2(\OO;\R^d)$ and $\sup_n \vert u_n \vert_n \leq c$. Then,
	 	\begin{align*}
	 		\int_\OO \Dfr_n(u_n(x)) \, \Dfr_n(v(x)) \dd x \to    \int_{\Omega}  \div u(x) \div v(x) \dd x.    
	 	\end{align*}
	 \end{corollary}

	Next, we investigate the limit of $(u_n, v_n)_n$. As $(\cdot, \cdot)_n$ can be interpreted as a nonlocal $H^1$-inner product,   one     can only expect convergence under strong   convergence    assumptions on $u_n$ and $v_n$. We will therefore   first    prove the result   in     the smooth case. Based on this, we can then compute the limit of $(u_n,v)_n$ for a weakly convergent sequence $u_n \wto u$ in $L^2(\OO;\R^d)$ and $v$ fixed, by resorting to mollification.     As the nonlocal operator considered in \cite{DM_bondbased} resembles the nonlocal inner product, the following lemma is reminiscent of   the results in \cite{DM_bondbased}.   

\begin{lemma}\label{lem:convergence smooth}
		Let $u_n, v_n \in C^2(\OO;\R^d)$ and $u,v \in   H^1_0(\Omega;\R^d)    $ be such that $\D u_n \to \D u$ and $\D v_n \to \D v$ strongly in $L^2(\OO;\R^{d \times d})$. Moreover, assume that  $\sup_n \Vert u_n \Vert_{W^{2,\infty}(\OO)^d} +  \Vert v_n \Vert_{W^{2,\infty}(\OO)^d} \leq c $.  Then, 
	\begin{align*}
		&\lim_{n \to \infty} \int_\OO	\int_\OO \rho_n (x'-x) \Dsc(u_n)(x', x)  \Dsc(v_n)(x',x)  \dd x' \dd x \\
		&=  d     \int_\OO \fint_{S^{d-1}} \left(s\cdot \D u(x) s \right) \left(s\cdot \D v(x) s \right)\dd\Hd^{d-1}(s)\dd x.
	\end{align*}
\end{lemma}
\begin{proof}
	Fix $x \in \OO$. Since the kernels localize around $x$, we split the integral into two parts. To this end let $R := \dist(x, \partial\OO)$ and $B(x,R)$ be the ball with radius $R$ centered at $x$. Then we have
	\begin{align}\label{eq:splitting domain of integration}
		&\int_\OO \rho_n(x'-x)  \Dsc(u_n)(x',x) \Dsc(v_n)(x',x) \dd x'  \notag \\
		&=   \left[ \int_{B(x,R)} + \int_{\OO\setminus B(x,R)} \right]     \rho_n(x'-x)\Dsc(u_n)(x', x) \Dsc(v_n)(x',x)  \dd x'. 
	\end{align}
	Note that, by Assumption~\ref{assumptions} on the kernels,  the second integral vanishes as $n\to \infty$, because
	\begin{align*}
		&\int_{\OO\setminus B(x,R)} \rho_n(x'-x) \Dsc(u_n)(x',x) \Dsc(v_n)(x',x) \dd x'  \\
		&\leq c(x) \Vert u_n \Vert_{L^\infty(\OO)^d} \Vert v_n \Vert_{L^\infty(\OO)^d}  \int_{\OO\setminus B(x,R)} \rho_n(x'-x) \dd x' \to 0. 
	\end{align*}
	We therefore restrict the analysis to   the     integral over $B(x,R)$, where we use a Taylor expansion   for      $\Dsc(u_n)$ and $\Dsc(v_n)$.   There exist  $z = z(x',x,u_n) \in B(x,R)$ and $z' = z'(x',x,v_n) \in B(x,R)$ such that     
	\begin{align*}%\label{eq:Taylor}
		\Dsc(u_n)(x',x) &=  \frac{x'-x}{|x'-x|} \cdot \D u_n(x) \frac{x'-x}{|x'-x|}  +   \frac12 (x' - x) \cdot \D^2 u_n(z) : \left( \frac{x'-x}{|x'-x|} \otimes \frac{x'-x}{|x'-x|} \right),     \\
		\Dsc(v_n)(x',x) &= \frac{x'-x}{|x'-x|} \cdot \D v_n(x) \frac{x'-x}{|x'-x|}    +\frac12 (x' - x) \cdot \D^2 v_n(z') : \left( \frac{x'-x}{|x'-x|} \otimes \frac{x'-x}{|x'-x|} \right),     
	\end{align*}
	  where $\D^2 u_n$  and  $\D^2 v_n$ are  the third order tensors of second derivatives.     
	  Using the above expressions,    by expanding the product $\Dsc(u_n)(x',x) \Dsc(v_n)(x',x)$ in the integral, we obtain four summands,   namely,    
	\begin{align*}
		\int_{B(x,R)} \rho_n(x'-x)\Dsc(u_n)(x', x) \Dsc(v_n)(x',x)  \dd x' = I_1 + I_2(u_n,v_n) + I_2(v_n,u_n) + I_3, 
	\end{align*}		
	where,   setting $\xi = x' - x$,    
	\begin{align}
		I_1 &=   \int_{B(0,R)}  \rho_n(\xi) \left( \frac{\xi}{|\xi|} \cdot \D u_n(x) \frac{\xi}{|\xi|}  \right) \left(  \frac{\xi}{|\xi|} \cdot \D v_n(x) \frac{\xi}{|\xi|}  \right)  \dd \xi, \label{eq:I1}  \\
		I_2(u_n,v_n) &=  \int_{B(0,R)}  \rho_n(\xi) \left( \frac{\xi}{|\xi|}\cdot \D u_n(x) \frac{\xi}{|\xi|}  \right)  \left(  \frac{\xi}{2}  \cdot \D^2 v_n(z') : \left( \frac{\xi}{|\xi|} \otimes \frac{\xi}{|\xi|}  \right) \right) \dd \xi, \label{eq:I2} \\
		I_3 &= 	\int_{B(0,R)} \rho_n(\xi)  \left(  \frac{\xi}{2}  \cdot \D^2 u_n(z) : \left( \frac{\xi}{|\xi|} \otimes \frac{\xi}{|\xi|}  \right) \right)  \left(  \frac{\xi}{2} \cdot \D^2 v_n(z') : \left( \frac{\xi}{|\xi|} \otimes \frac{\xi}{|\xi|}  \right) \right) \dd \xi.\label{eq:I3}
	\end{align}
 	To calculate the limit of $I_1$, note that $\xi \mapsto \left( \frac{\xi}{| \xi | } \cdot \D u_n(x) \frac{\xi}{| \xi | } \right) \left( \frac{\xi}{| \xi | } \cdot \D v_n(x) \frac{\xi}{| \xi | } \right)$ is positively homogeneous of degree $0$. Hence, as the kernels are radially symmetric by Assumption~\ref{assumptions},  we can write 
	\begin{align}\label{repi1}
		I_1 = \int_{B(0,R)}  \rho_n(\xi) \dd\xi \fint_{S^{d-1}} (s \cdot \D u_n(x) s) (s \cdot \D v_n(x) s) \dd\Hd^{d-1}(s),
	\end{align}
	see, e.g., \cite[Lemma~A.2]{kruzikQuasistaticElastoplasticityPeridynamics2018}. As $\D u_n \to \D u $ and $\D v_n \to \D v $ strongly in  $L^2(\OO;\R^{d \times d})$,  by using   Assumption~\ref{assumptions}~(ii),~(iii),     we therefore obtain
	\begin{align*}
		I_1 \to   d    \fint_{S^{d-1}} (s \cdot \D u(x) s)(s \cdot \D v(x) s) \dd\Hd^{d-1}(s),
	\end{align*}
	as $n \to \infty$. 

	The mixed terms $I_2(u_n,v_n)$ and $I_2(v_n,u_n)$ vanish as $n \to \infty$. Indeed,  
	\begin{align}\label{eq:product with grad vanishes}
		I_2(u_n,v_n) \leq c \Vert \D u_n \Vert_{L^\infty(\OO)^\dxd} \Vert \D^2 v_n \Vert_{L^\infty(\OO)^\dxdxd}    \int_{ B(0,R)   }\rho_n (\xi)  | \xi |   \dd\xi.
	\end{align}
	  Given    $0 < \epsi < R$,  we split the domain of integration and obtain
	\begin{align*}
 \int_{  B(0,R)   }\rho_n (\xi)  | \xi |   \dd\xi & =  \int_{B(0,\epsi)}\rho_n (\xi)  | \xi |   \dd\xi +  \int_{\epsi < \vert \xi \vert < R}\rho_n (\xi)  | \xi |   \dd\xi   \leq  \epsi d + R  \int_{\epsi < \vert \xi \vert < R} \rho_n (\vert \xi \vert) \dd\xi. 
	\end{align*}
The second term vanishes by Assumption~\ref{assumptions} and hence  by \eqref{eq:product with grad vanishes}  we conclude the claim as $\epsi> 0$   is     arbitrary   and the norms of $u_n$ and $v_n$ are  uniformly bounded. The analogous statement holds when exchanging the roles of $u_n$ and $v_n$, and hence $I_2(v_n,u_n) \to 0$ as $n \to \infty$,    as well.    	Moreover, the argument   above  also  implies that   $ I_3 \to 0 $	as $n \to \infty$.
	
	We thus conclude that 
	\begin{align*}
		\int_\OO \rho_n (x'-x) \Dsc(u_n)(x', x)  \Dsc(v_n)(x',x)  \dd x' \to   d     \fint_{S^{d-1}} \left(s\cdot \D u(x) s \right) \left(s\cdot \D v(x) s \right)\dd\Hd^{d-1}(s)
	\end{align*}
	pointwise in $x$ as $n \to \infty$. By the Dominated Convergence Theorem this also holds in $L^1(\OO;\R^d)$. Indeed, as $u_n$ and $v_n$ are smooth and their gradients are uniformly bounded, we can find   $M> 0$ such that
	\begin{align*}
		\left\vert (u_n(x') - u_n(x))\cdot	\frac{x'-x}{\vert x' - x\vert^2} \right\vert \leq M \quad \textup{ and }\quad 	\left\vert (v_n(x') - v_n(x))\cdot	\frac{x'-x}{\vert x' - x\vert^2} \right\vert \leq M, 
	\end{align*}
	and thus
	\begin{equation*}
		\int_{\OO} \rho_n(x'-x) \vert \Dsc(u_n)(x',x) \vert  \vert \Dsc(v_n)(x',x) \vert\dd x' \leq M^2 \int_{\OO} \rho_n(x' - x) \dd x' \leq   M^2 d.      \qedhere
	\end{equation*}
\end{proof}

  We now formulate two consequences of the above lemma.  

\begin{lemma}\label{lem:convergence time}
	Let $  w_n, w  \in L^2(0,T;H^1(\OO;\R^d))$ be such that $w_n(t) \in C^2(\OO;\R^d)$  for almost all $t \in [0,T]$. Moreover, assume that  $ w_n|_\Omega  \wto  w|_\Omega$ weakly in  $L^2(0,T;H^1(\Omega;\R^{d}))$ and that   $\sup_n \int_0^T \Vert w_n(t) \Vert^2_{W^{2,\infty}(\OO)^d} \dd t \leq c $.   Then,  we have %\RRR {\tt changed $A$ to $\Omega$ without making blue.}     
	\begin{align}\label{a new}
		&\liminf_{n \to \infty} \int_0^T \int_{\Omega  }	\int_{\Omega }	 \rho_n (x'-x) \Dsc(w_n)^2(x', x)    \dd x' \dd x  \dd t \notag\\
		&\geq   d    \int_0^T  \int_{\Omega  }	 \fint_{S^{d-1}} \left(s\cdot \D w(x) s \right)^2\dd\Hd^{d-1}(s)\dd x  \dd t. 
	\end{align}
\end{lemma}

\begin{proof}
The proof runs along the lines of the proof of Lemma~\ref{lem:convergence smooth} and we therefore only highlight the    adaptations.  Recalling the splitting \eqref{eq:splitting domain of integration}--\eqref{eq:I3}, the second integral  in  \eqref{eq:splitting domain of integration} and $I_2$, $I_3$ can be controlled uniformly in $x \in \OO$ in terms of $\eta_n  \Vert w_n(t) \Vert^2_{W^{2,\infty}(\OO)^d} $, where $(\eta_n)$ is a sequence with $\eta_n \to 0$ as $n \to \infty$. Thus, recalling the representation in \eqref{repi1} this shows
\begin{align*}
G  =  \liminf_{n \to \infty} \int_0^T \int_{\Omega  }	   \int_{B(0,R(x))}  \rho_n(\xi) \dd\xi \fint_{S^{d-1}} (s \cdot \D w_n(x) s)^2   \dd\Hd^{d-1}(s)   \dd x  \dd t,
\end{align*}
where $R(x)  = \dist(x,  \partial \OO 	)$ and $G$ denotes the left-hand side of \eqref{a new} for brevity. In view of 
 Assumption~\ref{assumptions}~(ii),~(iii), given any $\eps>0$, we get  
\begin{align*}
G \ge   \liminf_{n \to \infty}\int_0^T \int_\Omega   d(1-\eps ) \fint_{S^{d-1}}   (s \cdot \D w_n(x) s)^2   \dd\Hd^{d-1}(s)   \dd x  \dd t.
\end{align*}
As	$v \mapsto  \int_\Omega \fint_{S^{d-1}} (s \cdot \D v(x) s )^2 \dd \mathcal{H}^{d-1}(s) \dd x$ is  convex,  using  $  w_n|_\Omega  \wto  w|_\Omega $ weakly in $L^2(0,T;H^1(\Omega;\R^d))$  by assumption, we  obtain
	\begin{align*}
G \ge d(1-\eps) \int_0^T \int_\Omega \fint_{S^{d-1}} (s \cdot \D w (x) s)^2 \dd \Hd^{d-1}(s) \dd x \dd t .
	\end{align*}
As $\eps>0$ is arbitrary, this concludes the proof.  
\end{proof}

\begin{lemma}\label{lem:prod of nonloc strain}
	Let $u_n  \in \Sn(\OO;\R^d)$ and $u,v \in  H^1_0(\Omega;\R^d) $ be such that $u_n \to u$ strongly in $L^2(\OO;\R^d)$ and $\sup_n \vert u_n \vert_n < \infty$. Then,
	\begin{align*}
		\int_\OO \rho_n (x'-\cdot ) \Dsc(u_n)(x', \cdot)  \Dsc(v)(x',\cdot)  \dd x' \to    d     \fint_{S^{d-1}} \left(s\cdot \D u s \right) \left(s\cdot \D v s \right)\dd\Hd^{d-1}(s)
	\end{align*}
	strongly in $L^1(\OO;\R^d)$ as $n \to \infty$.
\end{lemma}

\begin{proof}
	\noindent\textit{Step 1.}   Assume momentarily that  $v \in C^\infty_c (\Omega;\R^d) $.  Let   $\eta \in C^\infty_c(\R^d)$, with $\Vert \eta \Vert_{L^1(\R^d)} = 1$ and set $\eta_k(x):= k^d \eta(kx)$. Moreover, set  $u_n^k :=  \eta_k  \ast u_n$ and $u^k :=  \eta_k   \ast u$,   where we implicitly extended $u_n$ and $u$ to $\R^d$ by $0$.    We have that $u_n^k \to u^k$ strongly in $L^2(  \OO; \R^d)$, because mollification is a compact operator, see, e.g., \cite[Corollary~4.28]{brezis}. Moreover, due to the mollification, also $u_n^k \to u_n$ strongly in $L^2(\OO;\R^d)$, and $u^k \to u$  strongly in $H^1( \OO;  \R^d)$.

	Then, we have that
	\begin{align*}
		&\left\Vert \int_\OO \rho_n (x'-\cdot) \Dsc(u_n)(x',\cdot) \Dsc(v)(x',\cdot)  \dd x' - \fint_{S^{d-1}} (s \cdot \D u(\blank) s ) (s \cdot \D v(\blank) s )\dd \Hd^{d-1}(s) \right\Vert_{L^1(\OO)} 
	\end{align*}
	 is controlled from above by the sum of three terms $I_1$, $I_2$, and $I_3$, defined as 
	\begin{align*}
		I_1 &= 	\left\Vert \int_\OO \rho_n (x'-\cdot) \Dsc(u_n)(x',\cdot)  \Dsc(v) (x',\cdot)\dd x' {-} \int_\OO \rho_n (x'-\cdot) \Dsc(u_n^k)(x',\cdot)  \Dsc(v)(x',\cdot)\dd x'  \right\Vert_{L^1(\OO)} \\
		I_2 &=	\left\Vert \int_\OO \rho_n (x'-\cdot) \Dsc(u_n^k)(x',\cdot)   \Dsc(v)(x',\cdot) \dd x'  {-} \fint_{S^{d-1}} \hspace{-0.1cm} (s \cdot \D u^k(\blank)  s)  (s \cdot \D v(\blank)  s) \dd \Hd^{d-1}(s) \right\Vert_{L^1(\OO)} \\ 
		I_3 &=	\left\Vert \fint_{S^{d-1}}   \hspace{-0.1cm} (s \cdot \D u^k(\blank)  s) (s \cdot \D v(\blank) s) \dd \Hd^{d-1}(s)  {-} \fint_{S^{d-1}}  \hspace{-0.1cm} (s \cdot \D u (\blank)s)  (s \cdot \D v (\blank)s ) \dd\Hd^{d-1}(s) \right\Vert_{L^1(\OO)}.
	\end{align*}
	Now, we let  first $n \to \infty$ and  then $ k \to \infty$ and show that each term   converges to $0$.     
	
	To check this for $I_1$, let us first note that, for all $w \in \Sn(\OO;\R^d)$, by   writing
	 $$\Dsc(w)(x',x) = w(x') \cdot \frac{x'-x}{ \vert x' - x\vert^2 } - w(x) \cdot \frac{x'-x}{  \vert x' - x\vert^2 } = -  w(x') \cdot \frac{x-x'}{  \vert x' - x\vert^2 } - w(x) \cdot \frac{x'-x}{  \vert x' - x\vert^2 }, $$ 
	 using a change of variable $x' \mapsto x$, the symmetry of $\rho_n$ and $\Dsc(v)$, and Fubini's Theorem,  we get    
	\begin{align*}%\label{eq:ibp}
		&\int_\OO \int_\OO \rho_n(x' - x)	\Dsc(w)(x', x) \Dsc(v) (x',x) \dd x' \dd x \notag \\
		&=  \int_\OO w(x) \cdot \left(   -    \int_\OO 2 \rho_n(x'-x) \frac{x'-x}{\vert x'-x \vert^2} \Dsc(v) (x',x) \dd x' \right) \dd x.	
	\end{align*}
	Setting $\Dsc_n^* \Dsc(v)(x) :=   -     \int_\OO 2 \rho_n(x'-x) \frac{x'-x}{\vert x'-x \vert^2} \Dsc(v) (x',x) \dd x'$,  H\"older's inequality  gives
	\begin{align*}
		I_1 &= \int_{ \OO } \left\vert \int_\OO \rho_n(x'-x) \Dsc(u_n - u_n^k)(x',x) \Dsc(v)(x',x) \dd x' \right\vert \dd x \\
		&\leq \int_{ \OO }  \vert(u_n - u_n^k)(x)\vert \, \vert \Dsc_n^*\Dsc(v)(x) \vert \dd x \leq \Vert u_n -u_n^k \Vert_{L^2(\OO)^d} \Vert  \Dsc_n^*\Dsc(v) \Vert_{L^2(\OO)^d}.
	\end{align*}
	Since $u_n \to u$, $u_n^k \to u^k$,   and $u^k \to u$     strongly in $L^2(\OO;\R^d)$ as $n\to \infty$, we are left with checking that $ \Vert  \Dsc_n^*\Dsc(v) \Vert_{L^2(\OO)^d}$ is uniformly bounded in $n$.
	  As $v$ is smooth and $\Dsc(v)(x',x) = \Dsc(v)(x,x')$, we obtain by Taylor expansion around $x$ and $x'$
	\begin{align*}
		2 \Dsc(v)(x',x) =  \frac{\xi}{| \xi |} \cdot \left( \D v(x') + \D v(x) \right) \frac{\xi}{| \xi |} + \frac12 \xi \cdot  \left( \D^2 v(z') + \D^2 v(z) \right) : \left( \frac{\xi}{| \xi |} \otimes \frac{\xi}{| \xi |} \right)
	\end{align*}   
	 for suitable points $z$ and $z'$, 	where we set $\xi := x' -x $ for brevity.  Accordingly, we split $ - \Dsc_n^* \Dsc(v)(x)  = J_1 + J_2$, with 
	\begin{align*}
		J_1 &:=  \int_\OO \rho_n(x' - x )  \frac{x'- x}{\vert x' -x \vert^2}  \frac{x' - x}{\vert x' - x\vert} \cdot \left( \D v(x') + \D v(x) \right) \frac{x' - x}{\vert x' - x\vert}\dd x', \\
		J_2 &:=	 \int_\OO \rho_n(x' - x )  \frac{x'- x}{\vert x' -x \vert^2}\frac{x' - x}{2} \cdot \left( \D^2 v(z') + \D^2 v(z) \right) : \left( \frac{x' - x}{\vert x' - x\vert}\otimes  \frac{x' - x}{\vert x' - x\vert} \right) \dd x'.
	\end{align*}
	To obtain an upper  bound for $J_2$, we use Assumption~\ref{assumptions}~(ii) to get  
$${\vert J_2 \vert  \le    \Vert \D^2 v \Vert_{L^\infty(\OO)^\dxdxd}  	\int_{\R^d}  \rho_n(\xi)   \dd \xi \leq d \Vert \D^2 v \Vert_{L^\infty(\OO)^\dxdxd}. } $$	
 	 To control $J_1$, we proceed analogously to \cite[Formula~(2.3)]{mengeshaLocalizationNonlocalGradients2015}.  For $x \in \OO$, we set $R(x) := \dist (x, \partial\OO)$ and  let  $\epsi \in (0, R(x))$ be arbitrary. Denoting by $\chi_M$ the characteristic function of the set $M$, we have for any $x' \in \OO$  
	 \begin{align*}%\label{eq:aux epsi 1}
- \chi_{[\epsi, \infty)}  (| x' -x |)  + 2 \chi_{[\epsi,R(x))}(| x' -x |)   	 	= \chi_{[\epsi, \infty)}(| x' -x |)    - 2 \chi_{[R(x), \infty)} (| x' -x |)  , 
	 \end{align*} 
	 and therefore also
	 \begin{align}\label{eq:J1 1}
	 	- &\chi_{[\epsi, \infty)} (| x' -x |)  \rho_n (x' - x )   \frac{x'-x}{|x'-x|^2 }  \frac{x'-x}{|x'-x| } \cdot (  \D v(x') + \D v(x) ) \frac{x'-x}{|x'-x| } \nonumber  \\
	 	&+ 2 \chi_{[\epsi,R(x))}(| x' -x |)  \rho_n (x' - x )   \frac{x'-x}{|x'-x|^2 }  \frac{x'-x}{|x'-x| } \cdot  \D v(x)\frac{x'-x}{|x'-x| }  \nonumber  \\
	 	=&	\chi_{[\epsi, \infty)}(| x' -x |)  \rho_n (x' - x )   \frac{x'-x}{|x'-x|^2 }  \frac{x'-x}{|x'-x| } \cdot \big(-(  \D v(x') - \D v(x)) \big) \frac{x'-x}{|x'-x| }  \nonumber \\
	 	&- 2 \chi_{[R(x), \infty)} (| x' -x |) \rho_n (x' - x )   \frac{x'-x}{|x'-x|^2 }  \frac{x'-x}{|x'-x| } \cdot  \D v(x)  \frac{x'-x}{|x'-x|^2 }. 
	 \end{align}
	  Since $  \xi \mapsto    \frac{\xi}{|\xi|^2 }  \frac{\xi}{|\xi| } \cdot  \D v(x)\frac{\xi}{|\xi| } $ is odd, integrating over a symmetric domain  around $x$  will give $0$. Therefore,   we have that 
	  \begin{align*}
	  	&2\int_\OO  \chi_{[\epsi,R(x))}(| x' -x |)  \rho_n (x' - x )   \frac{x'-x}{|x'-x|^2 }  \frac{x'-x}{|x'-x| } \cdot  \D v(x)\frac{x'-x}{|x'-x| } \dd x'  =  0. 
	  \end{align*}
%	  \\
%	  	&= 2 \int_\epsi^{R(x)} \tilde{\rho}_n(r) \int_{\partial B(0,r)}   \frac{x'-x}{|x'-x|^2 }  \frac{x'-x}{|x'-x| } \cdot  \D v(x)\frac{x'-x}{|x'-x| } \dd \Hd^{n-1}(x' ) \dd r =
	 In particular, integrating \eqref{eq:J1 1} with respect to $x'$ gives
	 \begin{align*}
	 	-& \int_\OO  \chi_{[\epsi, \infty)} (| x' -x |)  \rho_n (x' - x )   \frac{x'-x}{|x'-x|^2 }  \frac{x'-x}{|x'-x| } \cdot (  \D v(x') + \D v(x) ) \frac{x'-x}{|x'-x| } \dd x'  = H_1 + H_2  
	 \end{align*}
 with 	 
\begin{align*}
H_1 & = \int_\OO \chi_{[\epsi, \infty)}(| x' -x |)  \rho_n (x' - x )   \frac{x'-x}{|x'-x|^2 }  \frac{x'-x}{|x'-x| } \cdot \big(-(  \D v(x') - \D v(x)) \big) \frac{x'-x}{|x'-x| } \dd x',  \\
	 H_2 &  = - 2 \int_\OO  \chi_{[R(x), \infty)} (| x' -x |) \rho_n (x' - x )   \frac{x'-x}{|x'-x|^2 }  \frac{x'-x}{|x'-x| } \cdot  \D v(x)  \frac{x'-x}{|x'-x|^2 } \dd x'. 
	 \end{align*}
	 As $\D v$ is smooth, we can bound  
	 \begin{align*}
|H_1|	 	&\leq \int_\OO \chi_{[\epsi, \infty)}(| x' -x |)  \rho_n (x' - x )     \frac{1 }{|x'-x|}   \vert  \D v(x') - \D v(x)  \vert\dd x' \\
	 	&\leq \Vert \D^2 v \Vert_{L^\infty(\OO)^\dxdxd}   \int_{\R^d} \rho_n (\xi ) \dd \xi   \leq d \Vert \D^2 v \Vert_{L^\infty(\OO)^\dxdxd}. 
	 \end{align*}
	 We are left with estimating  $H_2$  
	 from above. To this end, pick $\xi(x) \in (\OO \setminus \supp \D v) \cap B(x, R(x))$. Then, $\D v(\xi(x)) = 0$ and by the choice of $\xi(x)$ also $\vert \xi(x) - x \vert \leq R(x)$. 
	 Therefore, we get
	 \begin{align*}
	 	 \big|\frac{1}{2}H_2\big|  	&\leq  \int_\OO  \chi_{[R(x), \infty)}(|x'-x|)  \rho_n (x' - x )   \frac{1}{|x'-x|}   \left| \D v(\xi(x)) - \D v(x) \right| \dd x' \\
	 	&\leq  \Vert  \D^2 v \Vert_{L^\infty(\OO)^\dxdxd}  \int_\OO  \chi_{[R(x), \infty)}(|x'-x|)  \rho_n (x' - x )    \frac{\vert \xi(x) - x \vert }{|x'-x| }  \dd x'  \leq d  \Vert  \D^2 v \Vert_{L^\infty(\OO)^\dxdxd} 
	 \end{align*}
	 since from the additional constraint $\chi_{[R(x), \infty)}(|x'-x|)$ on the integration domain, we have that $\vert \xi(x) - x \vert \leq R(x) \leq \vert x' - x\vert$. 
	 Hence, as $\epsi$ is arbitrary, we conclude that $\vert J_1 \vert \leq  3  d  \Vert  \D^2 v \Vert_{L^\infty(\OO)^\dxdxd}$. In particular,  the uniform $L^2$-bound of  $\Dsc_n^*\Dsc(v) $ follows.

	For showing that the term $I_2$ converges to $0$,  we can apply Lemma~\ref{lem:convergence smooth} as $u_n^k$ and $u^k$ are smooth.  Indeed, by Young's   Convolution    Inequality and Proposition~\ref{prop:poincare}, we have 
	\[ \sup_{n\in \N} \Vert u_n^k \Vert_{L^\infty(\OO)^d} \leq   \sup_{n\in \N} \Vert u_n \Vert_{L^2(\OO)^d} \Vert \eta_k   \Vert_{L^2(\R^d)} \leq c_k    \] 
	and analogously for $\D u_n^k$   and $\D^2 u_n^k$.  Moreover, this bound and Rellich's Theorem also  imply    that $\D u^k_n \to \D u^k $ strongly in $L^2(\OO; \R^\dxd)$.     

	  Eventually,     the third term $ I_3$  vanishes as $u^k \to u$ in $H^1(  \OO;  \R^d)$ by mollification. 
	
	\noindent\textit{Step 2.}   Let us now consider the case   $v \in H^1_0(\Omega;\R^d)$.     We approximate $v$ by a smooth sequence   $v^k \in C^\infty_c(\Omega;  \R^d)$,     i.e.,   we assume     $\Vert v - v^k \Vert_{ H^1_0(\Omega)^d} \to 0$, and  estimate 
	\begin{align}
		&\left\Vert \int_\OO \rho_n (x'-\cdot) \Dsc(u_n)(x', \cdot)  \Dsc(v)(x',\cdot)  \dd x' - \fint_{S^{d-1}} \left(s\cdot \D u(\blank) s \right) \left(s\cdot \D v(\blank) s \right)\dd\Hd^{d-1}(s) \right\Vert_{L^1(\OO)} \notag \\ 
		&\leq \int_\OO \int_\OO \rho_n (x' - x) \vert \Dsc(u_n)(x',x) \vert \, \vert \Dsc(v - v^k)(x', x)\vert  \dd x ' \dd x  \notag \\
		&+ \left\Vert \int_\OO \rho_n(x'-\cdot) \Dsc(u_n)(x',\cdot )  \Dsc(v^k)(x',\cdot)  \dd x' -   \fint_{S^{d-1}} (s \cdot \D  u  (\blank)  s ) (s \cdot \D v^k(\blank) s )\dd \Hd^{d-1}(s) \right\Vert_{L^1(\OO)} \notag \\
		&+ \int_\Omega  \fint_{S^{d-1}} | s \cdot \D u(x) s |  \, |s \cdot (\D v^k(x) - \D v(x) )s  | \dd \Hd^{d-1}(s) \dd x. \label{eq:step2}
	\end{align}
	We let first $n \to \infty$ and then $k \to \infty$. As  $v \in H^1_0(\Omega;\R^d)$,  we can apply Lemma~\ref{lem:bound Sn seminorm by H1},    H\"older's inequality, and the definition \eqref{eq:defintion nonlocal seminorm}      to  estimate the first term   on the right hand side of \eqref{eq:step2}     as
	\begin{align*}
		\int_\OO \int_\OO \rho_n (x' - x) \vert \Dsc(u_n)(x',x) \vert \, \vert \Dsc(v - v^k)(x', x)\vert  \dd x ' \dd x \leq \Vert v - v^k \Vert_{ H^1_0(\Omega)^d} \sup_{n\in \N}   \vert u_n \vert_n,     
	\end{align*}
 which tends to $0$ as $\sup_n \vert u_n \vert_n \leq c$ by assumption. 
	The second term    on the right hand side of \eqref{eq:step2}    vanishes by Step~1, whereas the third term     on the right hand side of \eqref{eq:step2}   converges to $0$ as $v^k$ converges to $v$ in $H^1(\Omega;\R^d)$.  
\end{proof}

\begin{remark}[Homogeneous Dirichlet boundary conditions]
	For the proof of the theorem above, it is instrumental to have homogeneous Dirichtlet boundary conditions, as otherwise one would only be able to establish an $L^1$-bound of $\Dsc_n^*\Dsc(v)$, see \cite[Proposition A.1]{neumann} and \cite[Proposition 2.12]{du}. 
\end{remark}

  As already mentioned, one     cannot expect   Lemma~\ref{lem:prod of nonloc strain}     to hold, if we replace $v$ by a sequence $v_n \wto v$.   However, one can still prove a lower semicontinuity result,   see, e.g.,      \cite[Lemma~A.1]{kruzikQuasistaticElastoplasticityPeridynamics2018}, \cite[Proof of Theorem~1.1]{mengeshaVariationalLimitClass2015}. 

\begin{lemma}[{\cite[Lemma~A.1]{kruzikQuasistaticElastoplasticityPeridynamics2018}}]\label{lem:moll lowers energy}
	Let $\eta_k$ be a standard   sequence of     mollifiers,  $u \in \Sn(\OO;\R^d)$ and   define   $u^k := \eta_k \ast u$. Then, for any $A \subset \subset \OO$ measurable, we have
	\begin{align*}
		\int_A \int_A \rho_n (x' - x) (\Dsc(u^k)(x',x))^2 \dd x' \dd x \leq \int_\OO \int_\OO \rho_n(x' - x) (\Dsc(u)(x', x))^2 \dd x' \dd x
	\end{align*} 
	   for all $k$ sufficiently large.    
\end{lemma} 
The proof directly follows from   the linearity of $u \mapsto \Dsc(u)$,  the convexity of $y \mapsto y^2$,  and Jensen's inequality.

%%%%%%%%%%%%%%%%%%%%%%%%%%%%%%%%%%%%%%%%%%%%%%%%%%%%%%%
	\section{Proof of Theorem~\ref{thm:main}}\label{sec:proof}
	%%%%%%%%%%%%%%%%%%%%%%%%%%%%%%%%%%%%%%%%%%%%%%%%%%%%%

	The proof is based on  the classical  evolutive $\Gamma$-convergence approach, inspired by \name{Sandier} and \name{Serfaty} \cite{sandierGammaconvergenceGradientFlows2004}, see also \cite{mielkeGamma}. This requires   to check the    $\liminf$-inequalities    \eqref{eq:li1}--\eqref{eq:li3},   which makes for the main part of the proof.
	
	\begin{lemma}[Liminf-inequalities]\label{lem:liminf inequalities}
		Let $u_n \in   H^1    (0,T;\Sn(\OO;\R^d))$ be a sequence of solutions to the nonlocal   viscoelastic   problem \eqref{eq:nonlocal problem},  and $u \in C([0,T];L^2(  \OO;  \R^d) )$ from the compactness result Proposition~\textup{\ref{prop:compactness}}.     Then,   extracting  a not relabeled subsequence with $u_n \to u$ in $C([0,T];L^2(\OO;\R^d))$,      for almost every $t \in [0,T]$, the following liminf-inequalities hold 
		\begin{align}
			&\liminf_{n \to \infty} E_n (u_n(t) ) \geq E(u(t)), \label{eq:liminf E}  \\
			&\liminf_{n \to \infty } \int_0^t D_n (\dot u_n(s)) \dd s  \geq \int_0^t D(\dot u(s) ) \dd s, \label{eq:liminf D} \\
			&\liminf_{n \to \infty } \int_0^t D_n^* \big(- \dd_n E_n (u_n(s)) \big) \dd s   \geq \int_0^t D^* \big(- \dd E(u(s)) \big) \dd s.  \label{eq:liminf D*}
		\end{align}
	\end{lemma}

	Note that the $\liminf$-inequality for the energy \eqref{eq:liminf E} has already been proved in \cite[Theorem~1.2]{mengeshaVariationalLimitClass2015} as part of the $\Gamma$-convergence result,  which we recalled in  Theorem~\ref{thm:gamma convergence}.  We  show  the two remaining liminf-inequalities \eqref{eq:liminf D} and \eqref{eq:liminf D*} in the subsequent Sections~\ref{subsec:D} and \ref{subsec:D star},   respectively.    
	
	%%%%%%%%%%%%%%%%%%%%%%%%%%%%%%%%%%%%%%%%%%%%%%%%%%%%%
	\subsection{Proof of {(\ref{eq:liminf D})} }\label{subsec:D}
	%%%%%%%%%%%%%%%%%%%%%%%%%%%%%%%%%%%%%%%%%%%%%%%%%%%%%
	 
	Let $(u_n)_n$ be a sequence of solutions to the nonlocal   viscoelastic   problem~\eqref{eq:nonlocal problem} and $u$ be the respective limit from the compactness result Proposition~\ref{prop:compactness}.   
	We resort to mollification and use Lemma~\ref{lem:moll lowers energy}.  As  $\Omega  \subset \subset \OO$,  for $k$ large enough  one has    
	\begin{align*}
		\int_{ \Omega} \int_{ \Omega} \rho_n(x' - x) (\Dsc( \eta^k \ast \dot{u}_n)(x',x))^2 \dd x ' \dd x \leq \vert \dot{u}_n \vert_n^2, 
	\end{align*}
	and hence also 
	\begin{align}\label{eq:Dn aux1}
		\liminf_{n\to\infty}   \int_0^T     \int_{ \Omega} \int_{ \Omega} \rho_n(x' - x) (\Dsc( \eta^k \ast \dot{u}_n)(x',x))^2 \dd x ' \dd x  \dd t \leq  \liminf_{n\to\infty}   \int_0^T     \vert \dot{u}_n \vert_n^2 \dd t. 
	\end{align}
   From the a priori bound \eqref{eq:unif bound dot u} and  Proposition~\ref{prop:poincare} we get that 
\begin{align}\label{eq:aux}
	\sup_{n\in \N} \int_0^T \Vert \dot u_n(t) \Vert_{L^2(\OO)^d}^2 \dd t \leq c.
\end{align}
Therefore,   one can find     a (non relabeled) subsequence $\dot{u}_n \wto v$ in $L^2(0,T;L^2(\OO;\R^d))$. By identifying weak limits, we see that $v = \dot u$. 
Moreover, \eqref{eq:aux} and   Young's Convolution Inequality implies    that $\sup_n \int_0^T  \Vert  ( \dot{u}_n (t) \ast \eta^k) \Vert_{W^{2,\infty}(\OO)^d}^2  \dd t  \leq c$. In particular, by a further weak compactness argument  this also implies that $\eta^k \ast \dot{u}_n \wto   \eta^k \ast  \dot{u}$ in $L^2(0,T;H^1(\Omega;\R^d))$.     Thus, we can apply Lemma~\ref{lem:convergence time} and obtain   
	\begin{align}\label{eq:Dn aux2}
		\liminf_{n\to\infty}   \int_0^T     \int_{ \Omega} \int_{ \Omega} \rho_n(x' - x) (\Dsc( \eta^k \ast \dot{u}_n)(x',x))^2 \dd x ' \dd x  \dd t \notag \\
		 \geq  d \int_0^T \int_{ \Omega} \fint_{S^{d-1}} (s \cdot   \D ( \eta^k \ast \dot u)(x)     s )^{2} \dd \Hd^{d-1}(s) \dd x   \dd t . 
	\end{align}
	  Combining \eqref{eq:Dn aux1} and \eqref{eq:Dn aux2}, we get    
	\begin{align*}
	 d 	\int_0^T  \int_\Omega  \fint_{S^{d-1}}  (s \cdot   \D ( \eta^k \ast  \dot{u})  (t,x)     s )^{2} \dd \Hd^{d-1}  (s)  \dd x \dd t  \leq \liminf_{n\to\infty}  \int_0^T  \vert  \dot{u}_n (t)  \vert_n^2 \dd t. 
	\end{align*}
	Now, we take the limit as $k\to \infty$ and   conclude as,  owing to the definition of $\vert \cdot \vert_n$ and Lemma~\ref{lem:vis potential forms},   the latter is exactly       inequality~\eqref{eq:liminf D}. \qed

	%%%%%%%%%%%%%%%%%%%%%%%%%%%%%%%%%%%%%%%%%%%%%%%%%%%%%
	\subsection{Proof of {(\ref{eq:liminf D*})} }\label{subsec:D star}
	%%%%%%%%%%%%%%%%%%%%%%%%%%%%%%%%%%%%%%%%%%%%%%%%%%%%%

  	The Fenchel Inequality gives for all $v \in H^1(\Omega;\R^d)$    
	\begin{align*}%\label{eq:D stern 0}
		D_n^* (- \dd_n E_n(u_n)) \geq        - \dd_n E_n(u_n)(v)     - D_n (v).
	\end{align*}
	The next lemma   shows that one can pass to the limit in the above   expression, which    is  the key ingredient to show \eqref{eq:liminf D*}. 
	\begin{lemma}[Convergence]\label{lem:lim}
		Let $u_n \in   H^1    (0,T;\Sn(\OO;\R^d))$ be a sequence of solutions to the nonlocal   viscoelastic    problem, $v \in H^1_0(\Omega;\R^d)$,  and $u \in C([0,T];L^2( \OO; \R^d) )$ from Proposition~\textup{\ref{prop:compactness}}.    Then,   extracting a non relabeled subsequence with $u_n \to u$ in $C([0,T];L^2(\OO;\R^d))$  we find,   as $n \to \infty$,    
		\begin{align}
			   - \dd_n E_n(u_n(t))(v)    &   \to      - \dd E(u(t))(v)     \quad \textup{for every $t \in [0,T]$}, \label{eq:D stern 1}\\
		D_n (v)  & \to D(v),\label{eq:D stern 2}
		\end{align}
  	where  $- \dd_n E(u_n)(v)$ and  $ - \dd E(u)(v) $ are as in \eqref{eq:derivative of energy} and \eqref{eq:gateaux dE}, respectively.       
	\end{lemma}
	\begin{proof}[Proof of Lemma~\ref{lem:lim}]
	 	The claim for $D_n$ \eqref{eq:D stern 2} follows from Lemma~\ref{lem:prod of nonloc strain}, when setting $u_n = v$   and recalling the various forms of $D$ in Lemma~\ref{lem:vis potential forms}.     
		
		For the    convergence of $- \dd_n E_n(u_n(t))$,     we expand the product in the integral of    \eqref{eq:derivative of energy}  as 
		\begin{align*}
			 - \dd_n E_n(u_n(t))(v) = &- \alpha \int_\OO \int_\OO \rho_n(x'- x) \Dsc(u_n(t))(x',x) \Dsc(v)(x',x) \dd x' \dd x \\
			 &- \left(\beta - \frac{\alpha}{d} \right) \int_\OO \int_\OO \Dfr_n(u_n(t))(x) \Dfr_n(v)(x) \dd x.  
		\end{align*}     Then we can employ Lemma~\ref{lem:prod of nonloc strain} and Corollary~\ref{cor:convergece of prod of div}   to pass to the limit.   Eventually, using  \eqref{eq:tensors intro},  Lemma~\ref{lem:prod  sDus},  and a computation analogous to the one of Lemma~\ref{lem:energy as tensor}   the proof is   concluded.    
	\end{proof}

	 We now proceed with the proof of   inequality   \eqref{eq:liminf D*} by first showing that $\liminf_{n \to \infty} D_n^* \big(- \dd_n E_n(u_n(t)) \big) \geq D^* \big(-\dd E(u(t)) \big)$    at almost every $t\in [0,T]$.   %This then allows us to conclude the proof of \eqref{eq:liminf D*} by employing the Fatou's Lemma by integrating in time. 

	   To this aim, let $\eta > 0$ be small but arbitrary. By definition of $D^*$  (see also \eqref{eq:defin D*}),  we can find $v \in H^1_0(\Omega;\R^d)$ such that 
	\begin{align*}
			D^* \big(-\dd E (u(t)) \big) \leq  - \dd E(u(t)) (v) - D(v) + \eta.
	\end{align*} 
Therefore,  by \eqref{eq:defin D*},  Lemma~\ref{lem:lim} implies 
	\begin{align*}
		\liminf_{n \to \infty } D_n^* \big(- \dd_n E_n (u(t)) \big) &\geq 	\liminf_{n \to \infty } \Big( - \dd_n E_n (u(t))(v) - D_n(v) \Big) \\
		&=  - \dd E(u(t))(v) - D(v) \geq 	D^* \big(-\dd E (u(t)) \big)  - \eta.
	\end{align*}
	Since $\eta$   is     arbitrary, we conclude   by     integrating over $t \in [0,T]$ and applying Fatou's Lemma.  \qed

	%%%%%%%%%%%%%%%%%%%%%%%%%%%%%%%%%%%%%%%%%%%%%%%%%%%%%
	\subsection{Conclusion of the proof Theorem \ref{thm:main}}\label{subsec:conclusion}
	%%%%%%%%%%%%%%%%%%%%%%%%%%%%%%%%%%%%%%%%%%%%%%%%%%%%%
	 	Let $u_n \in   H^1  (0,T;\Sn(\OO;\R^d))$ be a sequence of solutions to the nonlocal   viscoelastic   problem \eqref{eq:nonlocal problem},  and $u \in C([0,T];L^2(\OO;\R^d) ) $ from   Proposition~\textup{\ref{prop:compactness}}.

	From the well-preparedness assumption $E_n(u_n^0) \to E(u^0)$ and the energy-dissipation equality \eqref{eq:edi}, we obtain 
	\begin{align}\label{eq:conc 1}
		&E(u(0)) = \liminf_{n \to \infty} E_n(u_n(0))  \nonumber\\
		&\geq \liminf_{n \to \infty} E_n(u_n(t, \cdot))  + \liminf_{n \to \infty}  \int_0^t D_n(\dot u_n(s, \cdot)) \dd s   + \liminf_{n \to \infty} \int_0^t D_n^* (- \dd_n E_n(u_n(s, \cdot)) \dd s. 
	\end{align}
	Now, employing Lemma~\ref{lem:liminf inequalities}, we can further estimate the  right-hand side     of \eqref{eq:conc 1}    as
	\begin{align}\label{eq:conc 2}
		&\liminf_{n \to \infty} E_n(u_n(t, \cdot))  + \liminf_{n \to \infty}  \int_0^t D_n(\dot u_n(s, \cdot)) \dd s + \liminf_{n \to \infty} \int_0^t D_n^* (-\dd_n E_n( u_n(s, \cdot)) \dd s \nonumber \\
		&\geq E(u(t)) + \int_0^t D(\dot u(s, \cdot)) \dd s +  \int_0^t D^* (- \dd E(u(s, \cdot)) \dd s.
	\end{align}

	In order to conclude the proof, we need the classical chain rule for $E$. 
	
	\begin{lemma}[Chain rule]\label{lem:chain rule}
 	  	Let $\dot u \in L^2(0,T; H^1_0(\Omega;\R^d))$ and $\dd E(u) \in L^2(0,T;(H^1_0(\Omega;\R^d))^*)$. Then, the chain rule holds, i.e.,  
		\begin{align*}
			\frac{\dd }{ \dd t  } E(u(t)) =  \dd  E(u(t))( \dot u (t) ) \qquad \textup{for almost every } t \in [0,T]. 
		\end{align*}
	\end{lemma}

		 We now check the regularity  assumptions of the lemma above. 
		As $E(u(0)) < \infty$ and $D^* \geq 0$, we obtain from    \eqref{eq:conc 1}--\eqref{eq:conc 2} 
		\begin{align*}
			\int_0^t \int_\Omega \Dd \epsi(\dot u(t,x)) : \epsi (\dot u(t,x)) \dd x \dd t < \infty,
		\end{align*}
		which is nothing else but $\Vert \dot u\Vert_{L^2(0,T;H^1_0(\Omega;\R^d))} < \infty$ by  Lemma~\ref{lem:korn type}.  This, along with the fact that $u(0) \in H^1_0(\Omega;\R^d)$, also shows that  $u \in L^2(0,T;H^1_0(\Omega;\R^d))$. 

		Next, we prove that $\dd E(u) \in L^2(0,T;  (H^{1}_0(\Omega;\R^d))^*    )$. First note that from   \eqref{eq:gateaux dE}     we have
		\begin{align*}
			  \vert  \dd E(u)(v) \vert     \leq c \int_\Omega \epsi(u(x)) : \epsi(v(x)) \dd x   \leq c \Vert u\Vert_{H^1_0(\Omega)^d} \Vert v \Vert_{H^1_0(\Omega)^d} 
		\end{align*}
		for some constant $c > 0$. Therefore, we arrive at
		\begin{align*}
			\int_0^t \Vert \dd  E(u(s)) \Vert^2_{  (H^{1}_0(\Omega;\R^d))^*    } \dd s & = \int_0^t \left( \sup_{v \in H^1_0(\Omega;\R^d) : \Vert v \Vert_{H^1_0(\Omega)^d}\leq 1}  \vert    \dd  E(u)(v)     \vert \right)^2  \dd s  \nonumber\\
			&\leq c \int_0^t \left(  \sup_{v \in H^1_0(\Omega;\R^d) : \Vert v \Vert_{H^1_0(\Omega)^d}\leq 1}\Vert u \Vert_{H^1_0(\Omega)^d} \Vert v \Vert_{H^1_0(\Omega)^d}   \right)^2 \dd s 
		\nonumber \\ &	\leq c \int_0^t   \Vert u(t) \Vert_{H^1_0(\Omega)^d}^2 \dd s  
			 \le  c \Vert u \Vert^2_{L^2(0,T;H^1_0(\Omega;\R^d))}.  %\label{eq:somebound}
		\end{align*}
	 In particular, the assumptions of Lemma~\ref{lem:chain rule} are fulfilled.    Now,  from     \eqref{eq:conc 1}--\eqref{eq:conc 2}, and the chain rule    of     Lemma~\ref{lem:chain rule}, we obtain
	\begin{align*}
		0 &\geq E(u(t)) - E(u(0)) + \int_0^t D(\dot u(s)) \dd s +  \int_0^t D^* \big(- \dd E(u(s)) \big) \dd s \\
			&= \int_0^t \Big(    \dd E(u(s)) ( \dot u(s) )      + D(\dot u(s)) + D^* \big(  - \dd E (u(s)) \big)  \Big) \dd s.
	\end{align*}
Since the integrand	$     \dd E(u)(\dot u )     + D(\dot u) + D^* (-\dd E ( u)) $ is nonnegative due to the definition of the Fenchel-Legendre transform, we conclude     that    
	\begin{align*}%\label{eq:something}
		  	 \dd  E(u(t)) (\dot u(t) )       + D(\dot u(t)) + D^* \big( -\dd E ( u(t)) \big) = 0 	\qquad  \textup{ for almost every  }  t \in [0,T]. 
	\end{align*}

	Therefore, $u$ is a   (weak)    solution of \eqref{eq:pde}.

	Uniqueness follows  by standard arguments  using the linearity of the limit problem. 
%	Indeed, let $u$ and  $w$ be two solutions. Then  the energy dissipation equality for the limiting problem is satisfied for  $u-w$ as well and hence  \eqref{eq:bound 2} and standard arguments imply that $u = w$. 
	In particular,   all extractions of subsequences as $n \to \infty$ were unnecessary and      the whole sequence $u_n$ converges to $u$. 	 \qed
	
	% chain rule and solution => d/dt E < 0. da E(0) = 0 => E(u-w) = 0 für alle Zeiten und damit auch | u - w | = 0 für alle Zeiten. 
	
	From the above proof, we readily get the following   additional information.    
	
	\begin{corollary}[Energy  convergence]\label{cor: en conv}
		  Let $u_n$ be a sequence of solutions to the nonlocal   viscoelastic   problem~\eqref{eq:nonlocal problem} and $u$ be the solution to the local   viscoelastic   problem~\eqref{eq:pde} with $u_n \to u$ as in the setting of Theorem~\textup{\ref{thm:main}}.    Then, we have
		\begin{align*}
			E_n(u_n(t)) \to E(u(t)) \quad \textup{for almost all $t \in [0,T]$}. 
		\end{align*}
	\end{corollary}

			\section*{Acknowledgments} 
This work was supported by the DFG project FR 4083/3-1, by the
Deutsche Forschungsgemeinschaft (DFG, German Research Foundation)
under Germany's Excellence Strategy EXC 2044 -390685587, Mathematics
M\"unster: Dynamics--Geometry--Structure, by the Austrian
Science Fund (FWF) projects \href{https://doi.org/10.55776/F65}{10.55776/F65}, \href{https://doi.org/10.55776/I4354}{10.55776/I4354}, \href{https://doi.org/10.55776/I5149}{10.55776/I5149}, and \href{https://doi.org/10.55776/P32788}{10.55776/P32788}. 
	
	\appendix

	%%%%%%%%%%%%%%%%%%%%%%%%%%%%%%%%%%%%%%%%%%%%%%%%
	\section{Calculations}\label{sec:appendix calculations}
	%%%%%%%%%%%%%%%%%%%%%%%%%%%%%%%%%%%%%%%%%%%%%%%%
	
	 In this section, we collect some elementary calculations.

	\begin{lemma}[Integrals over the unit sphere]\label{lem:surface intergrals}
		\begin{align*}
			&\fint_{S^{d-1}} s_i^n s_j \dd\Hd^{d-1}(s) = 0\quad &&\forall i \neq j, i,j \in \{ 1, \dots, d\}, \, \forall n \in \N, \\
			&\fint_{S^{d-1}} s_i^2 \dd\Hd^{d-1}(s) = \frac1d \quad && \forall i \in \{1,\dots, d \}, \\
			&\fint_{S^{d-1}} s_i^4 \dd\Hd^{d-1}(s) = \frac{3}{d(d+2)}  \quad &&\forall i \in \{1,\dots, d \}, \\
			&\fint_{S^{d-1}} s_i^2 s_j^2 d\Hd^{d-1}(s) = \frac{1}{d(d+2)}  \quad &&\forall i \neq j, i,j \in \{ 1, \dots, d\}.
		\end{align*}
	\end{lemma}
	
	A proof can be found in   \cite[Theorem~A.1]{mengeshaNONLOCALKORNTYPECHARACTERIZATION2012},     see also \cite{baker}.
	
	\begin{lem}\label{lem:div}
		 For $u \in H^1(\Omega;\R^d)$ we have \[ \fint_{S^{d-1}} s \cdot \D u(x) s \dd\Hd^{d-1}(s) = \frac1d \div u(x)\qquad  \textup{ for a.e. $x \in \Omega$}. \]
	\end{lem}

	\begin{proof}
	  A simple computation shows     
		\begin{align*}
			\fint_{S^{d-1}} s \cdot \D u(x) s \dd\Hd^{d-1}(s)  &= \sum_{ij} \partial_i u_j (x) \fint_{S^{d-1}} s_i s_j  \dd\Hd^{d-1}(s) \\
			&= \sum_i \partial_i u_i (x) \frac1d = \frac1d \div u(x),
		\end{align*}
		where we used Lemma~\ref{lem:surface intergrals}.
	\end{proof}

	We can also give an expression for the product of quadratic forms on the sphere.
	
	\begin{lemma}\label{lem:prod  sDus}
			Let $u, v \in H^1(\Omega;\R^d)$. Then we have for almost every $x \in \Omega$
			\begin{align*}%\label{eq:int of product of quadratic forms}
				&\fint_{S^{d-1}} (s \cdot \D u(x) s) (s \cdot \D v(x) s) \dd\Hd^{d-1}(s) \\
				&=\frac{2}{d(d+2)} \epsi(u(x)): \epsi(v(x))+ \frac{1}{d(d+2)} (\div u(x) ) (\div v(x)).
			\end{align*}
	\end{lemma}
	\begin{proof}
		As for any matrix $M$ we have $x\cdot Mx = x \cdot \frac{M + M^\top}{2} x$, we can replace $\nabla u$ by the symmetric gradient $\epsi(u)$ in the following calculations. 
		For brevity we set $A := \epsi(u)$ and $B := \epsi(v)$. Then, we obtain by Lemma~\ref{lem:surface intergrals}
			\begin{align*}
			&\fint_{S^{d-1}} (s \cdot \D u s) (s \cdot \D v s) \dd \Hd^{d-1}(s) = \fint_{S^{d-1}} (s\cdot A s)(s \cdot Bs) \dd \Hd^{d-1}(s) \\
			&= \sum_{ijkl} A_{ij}  B_{kl}  \fint_{S^{d-1}} s_i s_j s_k s_l \dd\Hd^{d-1}(s) \\
			&=   \frac{3}{d(d+2)}\sum_i A_{ii} B_{ii} +  \frac{1}{d(d+2)}\sum_{i\neq k} A_{ii}B_{kk}+  \frac{1}{d(d+2)}\sum_{i\neq j} A_{ij} B_{ij}+ \frac{1}{d(d+2)}\sum_{i\neq j} A_{ij} B_{ji}.     
		\end{align*}
		Next, note that, since $\partial_i u_i = \epsi(u)_{ii}$, we have
		\begin{align*}
			\sum_{i} A_{ii} B_{ii} + \sum_{i \neq k} A_{ii} B_{kk} = \sum_{i,k} A_{ii} B_{kk}  = \left(  \sum_i \partial_i u_i \right)  \left(  \sum_k \partial_k v_k \right) =	\div u \,  \div v.
		\end{align*}
		Eventually, by virtue of the symmetry of $A$ and $B$, we have 
		\begin{align*}
			2 	\sum_ i A_{ii} B_{ii} +  \sum_{i \neq j} A_{ij} B_{ij} + \sum_{i\neq j} A_{ij} B_{ji} &=  2 \sum_i \sum_j A_{ij} B_{ji} = 2 \sum_i (AB)_{ii}   \\
			&= 2 \tr(AB^\top) =  2 A:B = 2  \epsi(u) : \epsi(v).
		\end{align*}
		Putting this together yields the claim.
	\end{proof}

	Let us conclude by proving the claim on the equivalent forms of the energy in Theorem~\ref{thm:gamma convergence}.  We present a full argument for the sake of completeness, but one could also refer to \cite[Proposition~4]{mengeshaNonlocalConstrainedValue2014}.
	
	\begin{lemma}\label{lem:energy rewritten}
		Let $u,v \in H^1(\Omega;\R^d)$. Then
		\begin{align*}
			 &\int_{\Omega} \fint_{S^{d-1}} \left(s \cdot \D u (x) s - \frac1d \div u(x) \right) \left(s \cdot \D v (x) s - \frac1d \div v (x) \right) \dd \Hd^{d-1}(s)  \dd x  \\
			 &=  \int_{\Omega} \frac{2}{d(d+2)} \epsi(u(x)) : \epsi(v(x)) +  \left(\frac{1}{d(d+2)} - \frac{1}{d^2} \right) (\div u(x) ) (\div v(x)) \dd x.
		\end{align*}
	In particular, we have
		\begin{align*}
		E (u) & =    \frac{\mu}{2}     \int_\Omega \vert \epsi(u(x)) \vert^2 \dd x	+   \frac{\lambda}{2}    \int_\Omega \div^2 u(x) \dd x \\
		& =  \frac{\beta}{2} \int_\Omega \div^2 u(x) \dd x  + \frac{\alpha}{2}  d     \int_\Omega \fint_{S^{d-1}} \left(s \cdot \D u (x) \cdot s - \frac1d \div u(x) \right)^2 \dd \Hd^{d-1}(s) \dd x,
	\end{align*}
	where $\mu = \frac{2\alpha}{d+2}$ and   $\lambda =\beta - \frac{2\alpha}{d (d+2)}$.     
	\end{lemma}

	\begin{proof}
		By Lemma~\ref{lem:div} we have
		\begin{align*}
			\int_\Omega \fint_{S^{d-1}} (s \cdot \D u (x) s )  \frac{1}{d} \div v (x)   \dd\Hd^{d-1}(s) \dd x  = \frac{1}{d^2} \int_\Omega \div u(x) \, \div v(x) \dd x
		\end{align*}
		 and similarly for $u$ and $v$ interchanged. We can hence conclude 
		\begin{align*}
			 &\int_{\Omega} \fint_{S^{d-1}} \left(s \cdot \D u (x) s - \frac1d \div u(x) \right) \left(s \cdot \D v (x) s - \frac1d \div v (x) \right) \dd\Hd^{d-1}(s) \dd x \\
			 &=	 \int_{\Omega} \fint_{S^{d-1}} \left(s \cdot \D u (x) s \right)  \left(s \cdot \D v (x) s \right) \dd\Hd^{d-1}(s) \dd x - \frac{1}{d^2} \int_\Omega  \div u(x) \, \div v(x) \dd x. 
		\end{align*}
		Employing Lemma~\ref{lem:prod  sDus}, we arrive at
		\begin{align*}
			&\int_{\Omega} \fint_{S^{d-1}} \left(s \cdot \D u (x) s \right)  \left(s \cdot \D v (x) s \right) \dd\Hd^{d-1}(s) \dd x - \frac{1}{d^2} \int_\Omega  \div u(x) \, \div v(x) \dd x \\
			&= \int_\Omega \left( \frac{2}{d(d+2)} \epsi(u): \epsi(v) + \frac{1}{d(d+2)} \div u(x ) \,  \div v(x) - \frac{1}{d^2} \div u(x) \,  \div v(x) \right)  \dd x
		\end{align*}
		 concluding the  first part of the  claim. The representation of the energy follows by   \eqref{eq:gamma limit of energy}.  		
	\end{proof}

\begin{lemma}[Elastic energy]\label{lem:energy as tensor}
	With $\C$ defined as in \eqref{eq:tensors intro}, we have 
	\begin{align*}%\label{eq:e with C}
		 \frac{\lambda}{2} \int_\Omega \div u(x)^2 \dd x + \frac{\mu}{2} \int_\Omega \vert \epsi(u(x)) \vert^2 \dd x =    \frac12     \int_\Omega \C \epsi(u(x)) : \epsi (u(x)) \dd x.
	\end{align*}
\end{lemma}
  \begin{proof}
	Set for brevity $\epsi_{ij} := \epsi(u)_{ij}$.  A  direct calculation shows that
	\begin{align*}
		\C \epsi : \epsi &=  \sum_{ijkl}   \C_{ijkl} \epsi_{ij} \epsi_{kl}  =  \sum_{ijkl}    \big( \mu   \delta_{ik} \delta_{jl} \epsi_{ij} \epsi_{kl} + \lambda   \delta_{ij} \delta_{kl} \epsi_{ij} \epsi_{kl} \big)\\
		&=   \sum_{ij}   \mu \epsi_{ij} \epsi_{ij} +  \sum_{ik}   \lambda  \epsi_{ii} \epsi_{kk}   = \mu \vert \epsi\vert^2 + \lambda (\div \epsi)^2. \qedhere
	\end{align*}
\end{proof}

	%BEGIN_FOLD Bibliography

	%END_FOLD
\end{document}